\documentclass[a4paper,11pt]{amsart}
\usepackage{graphicx,color}
\usepackage{amsmath, amsthm, amssymb, amsfonts}
\usepackage{amscd}
\usepackage[all]{xy}

\newtheoremstyle{Style}
{.5em}{.5em}
{\it}
{}
{\sc}
{\ {\bf---}}
{ }
{}

\newtheoremstyle{StyleRemarque}
{\topsep}{.5em}
{\it}
{}
{\slshape}
{.}
{ }
{}

\theoremstyle{Style}
\newtheorem{defn}{Definition}[section]
\newtheorem{lem}[defn]{Lemma}
\newtheorem{prop}[defn]{Proposition}
\newtheorem{thm}[defn]{Theorem}
\newtheorem{cor}[defn]{Corollary}
\newtheorem*{thm*}{Theorem}

\theoremstyle{StyleRemarque}
\newtheorem{rem}{Remark}
\newtheorem{ex}{Example}

\newcommand{\norm}[1]{\left\Vert#1\right\Vert}
\newcommand{\abs}[1]{\left\vert#1\right\vert}
\newcommand{\derivzero}[1]{\frac{d}{d#1}\Big\vert_{#1=0}}

\newcommand{\set}[1]{\left\{#1\right\}}
\newcommand{\Forall}[2]{\forall \, #1 \in #2, \:}
\newcommand{\Exists}[2]{\exists \, #1 \in #2, \:}
\newcommand{\diagcommut}[8]{
$$\begin{CD}
 #1@>{#2}>>#3  \\
 @VV{#4}V   @VV{#5}V \\
 #6@>{#7}>>#8 \\
\end{CD}$$}
\newcommand{\restric}[1]{\vert_{#1}}
\newcommand{\appli}[3]{#1 : \, #2 \To #3}

\newcommand{\Rl}{\mathbb R}
\newcommand{\Cx}{\mathbb C}
\newcommand{\Ir}{\mathbb Z}
\newcommand{\Nl}{\mathbb N}
\newcommand{\Sph}{\mathbb S}

\newcommand{\To}{\longrightarrow}
\newcommand{\U}{\mathcal{U}}

\newcommand{\A}{\mathcal{A}}

\newcommand{\OO}{\mathcal O}

\newcommand{\I}{\mathcal I}

\DeclareMathOperator{\Tr}{Tr}
\DeclareMathOperator{\Ric}{Ric}

\DeclareMathOperator{\Sect}{Sect}

\DeclareMathOperator{\vol}{vol}

\DeclareMathOperator{\inj}{inj}

\DeclareMathOperator{\Id}{Id}
\DeclareMathOperator{\grad}{grad}

\DeclareMathOperator{\id}{id}
\DeclareMathOperator{\Rm}{Rm}
\DeclareMathOperator{\Exp}{Exp}

\begin{document}


\title[On asymptotically flat manifolds with non maximal volume growth.]{\bfseries{On some asymptotically flat manifolds with non maximal volume growth.}}
\author{Vincent Minerbe}
\date{\today}

\maketitle

\begin{abstract}
Asymptotically flat manifolds with Euclidean volume growth are known to be ALE. In
this paper, we consider a class of asymptotically flat manifolds with slower volume growth
and prove that their asymptotic geometry is that of a fibration over an ALE manifold. In particular, we show that gravitational instantons with cubic volume growth are ALF.
\end{abstract}


\section*{Introduction.}

The aim of this paper is to understand the geometry at infinity of some asymptotically flat manifolds, that is complete noncompact Riemannian 
manifolds $(M,g)$ whose curvature tensor $\Rm_g$ fastly goes to zero at infinity:
$$
\abs{\Rm}_g = \OO(r^{-2-\epsilon})
$$
where $r$ is the distance function to some point and $\epsilon$ is a positive number. Such manifolds are known to have finite topological 
type \cite{A2}: there is a compact subset $K$ of $M$ such that $M \backslash K$ has the \emph{topology} of $\partial K \times \Rl_+^*$. In contrast, any (connected) manifold carries a complete metric with quadratic curvature decay ($\abs{\Rm}_g = \OO(r^{-2}$), see \cite{LS}). 

In this paper, we are interested in the \emph{geometry} at infinity of some asymptotically flat manifolds: we want to understand the metric $g$ in the unbounded part 
$M \backslash K$. 

A basic fact is that an asymptotically flat manifold has at most Euclidean volume growth: there is a constant $B$ such that
$$
\Forall{x}{M} \forall \, t \geq 1, \, \vol B(x,t) \leq B t^n.
$$
A fundamental geometric result was proved in \cite{BKN}: an asymptotically flat manifold with maximal volume growth, that is 
$$
\Forall{x}{M} \forall \, t \geq 1, \, \vol B(x,t) \geq A t^n,
$$
is indeed Asymptotically Locally Euclidean (``ALE'' for short): there is a compact set $K$ in $M$, a ball $B$ in $\Rl^n$, a finite subgroup $G$ of $O(n)$ and a diffeomorphism $\phi$ between $\Rl^n \backslash B$ and $M \backslash K$ such that $\phi^* g$ tends to the standard metric $g_{\Rl^n}$ at infinity. It is also proved in \cite{BKN} that a complete Ricci flat manifold with maximal volume growth and curvature in $L^{\frac{n}{2}}(dvol)$ is ALE. 

In the paper \cite{BKN}, S. Bando, A. Kasue and H. Nakajima raise the following natural question: can one understand the geometry at infinity of asymptotically flat manifolds whose volume growth is \emph{not} maximal ?

Some motivation comes from theoretical physics. Quantum gravity and string theory make use of the so called ``gravitational instantons'': from a mathematical point of view, a gravitational instanton is a complete noncompact hyperk\"ahler $4$-manifold with decaying curvature at infinity. In dimension four, ``hyperk\"ahler'' means Ricci flat and K\"ahler. Note that in the initial definition, by S. Hawking \cite{Haw}, gravitational instantons were not supposed to be K\"ahler, including for instance the Riemannian Schwarzschild metric. The curvature should satisfy a ``finite action'' assumption: typically, we want the curvature tensor to be in $L^2$. These manifolds have recently raised a lot of interest, both from mathematicians (for instance, \cite{HHM}, \cite{EJ}) and physicists. Works by \cite{BKN} and \cite{K1}, \cite{K2} enable to classify the geometry of gravitational instantons with maximal volume growth: these are ALE manifolds, obtained as resolutions of quotients of $\Cx^2$ by a finite subgroup of $SU(2)$. In case the volume growth is not maximal, S. Cherkis and A. Kapustin have made a classification conjecture, inspired from string theory (\cite{EJ}). Essentially, gravitational instantons should be asymptotic to fibrations over a Euclidean base and the fibers should be circles (``ALF'' case, for ``Asymptotically Locally Flat''), tori (``ALG'' case, because E,F,... G) or compact orientable flat $3$-manifolds (there are $6$ possibilities, this is the `` ALH '' case).  

Let us state our main theorem. Here and in the sequel, we will denote by $r$ the distance to some fixed point $o$, without mentionning it. We will 
also use the measure $d\mu = \frac{r^n}{\vol B(o,r)} dvol$. It was shown in \cite{Min} that this measure has interesting properties on 
manifolds with nonnegative Ricci curvature. Note that in maximal volume growth, it is equivalent to the Riemannian measure $dvol$. 

\begin{thm}
Let $(M^4,g)$ be a connected complete hyperk\"ahler manifold with curvature in $L^2(d\mu)$ and whose volume growth obeys
$$
\forall \, x \in M, \, \forall \, t \geq 1, \, A t^\nu \leq \vol B(x,t) \leq B t^\nu
$$
with $0 < A \leq B$ and $3 \leq \nu < 4$. Then $\nu=3$ and there is a compact set $K$ in $M$, a ball $B$ in $\Rl^3$, a finite subgroup $G$ of $O(3)$ and a circle fibration $\appli{\pi}{M \backslash K}{(\Rl^{3} \backslash B)/ G}$. Moreover, the metric $g$ can be written 
$$
g = \pi^*\tilde{g} + \eta^2 + \mathcal{O}(r^{-2}),
$$  
where $\eta^2$ measures the projection along the fibers and $\tilde{g}$ is ALE on $\Rl^3$, with 
$$
\tilde{g} = g_{\Rl^3} + \OO(r^{-\tau}) \quad \text{for every } \, \tau<1.
$$
\end{thm}

Up to finite covering, the topology at infinity (i.e. modulo a compact set) is therefore either that of $\Rl^3 \times \Sph^1$ (trivial fibration over $\Rl^3$) or that of $\Rl^4$ (Hopf fibration). Examples involving the Hopf fibration are provided by S. Hawkings' multi-Taub-NUT metrics, which we will describe in the text.

Our assumption on the curvature may seem a bit strange at first glance. It should be noticed it is a priori weaker than $\Rm = \OO(r^{-2-\epsilon})$. So we are in the realm of asymptotically flat manifolds. In the appendix, we will show that indeed, our integral assumption implies stronger estimates. Under our hypotheses, a little analysis provides $\nabla^k \Rm = \OO(r^{-3-k})$, for any $k$ in $\Nl$!   

Our volume growth assumption is uniform: the constants $A$ and $B$ are assumed to hold at any point $x$. This is not anecdotic. By looking at flat examples, we will see the importance of this uniformity. This feature is not present in the maximal volume growth case, where the uniform estimate    
$$
\exists \, A, \, B \in \Rl_+^*, \, \forall \, x \in M, \, \forall \, t \geq 1, \, A t^n \leq \vol B(x,t) \leq B t^n
$$
is equivalent to
$$
\exists \, A, \, B \in \Rl_+^*, \exists \, x \in M, \, \forall \, t \geq 1, \, A t^n \leq \vol B(x,t) \leq B t^n.
$$
  
Our result is related to a theorem by A. Petrunin and W. Tuschmann \cite{PT}. They have shown that if an asymptotically flat $4$-manifold has a simply connected end, then this end admits a tangent cone at infinity that is isometric to $\Rl^4$, $\Rl^3$ or $\Rl \times \Rl_+$. The $\Rl^4$ case is the ALE case and the case $\Rl \times \Rl_+$ conjecturally never occurs.  So we are left with $\Rl^3$, and this is consistent with our theorem (which does not require simply-connectedness at infinity, by the way). In greater dimension $n$, \cite{PT} says the tangent cone at infinity is $\Rl^n$ as soon as the end is simply connected. Indeed, our theorem admits a version (see \ref{thmgl} below) in any dimension and without the hyperk\"ahler assumption; it requires estimates on the covariant derivatives of the curvature, \emph{but also} an unpleasant assumption on the holonomy of short loops at infinity (it should be sufficiently close to the identity). Under these strong assumptions, a volume growth comparable to that of $\Rl^{n-1}$ in dimension $n \geq 5$ implies the existence of a necessarily trivial circle fibration over the complement of a ball in $\Rl^{n-1}$, which forbids simply connected ends. In this point of view, the fact that dimension $4$ is special in \cite{PT} comes from the existence of the Hopf fibration over $\Sph^2$.    

The idea of the proof is purely Riemannian. The point is the geometry at infinity collapses, the injectivity radius remains bounded while the curvature gets very small, so Cheeger-Fukaya-Gromov theory \cite{CG},\cite{CFG} applies. The fibers of the circle fibration will come from suitable regularizations of short loops based at each point. 

The structure of this paper is the following. 

In a first section, we will consider examples, with three goals: first, we want to explain our volume growth assumption through the study of flat manifolds; second, these flat examples will also provide some ideas about the techniques we will develop later; third, we will describe the Taub-NUT metric in detail, so as to provide the reader with a concrete example to think of.   

In a second section, we will try to analyse some relations between three Riemannian notions: curvature, injectivity radius, volume growth. We will 
introduce the ``fundamental pseudo-group''. This object, due to M. Gromov \cite{GLP}, encodes the Riemannian geometry at a fixed scale. It is our basic tool and its study will explain for instance the volume growth self-improvement phenomenon in our theorem (from $3\leq \nu < 4$ to $\nu=3$).  

In the third section, we completely describe the fundamental pseudo-group  at a convenient scale, for gravitational instantons. This enables us to build the fibration at infinity, locally first, and then globally, by a gluing technique. Then we make a number of estimates to obtain the description of the geometry at infinity that we announced in the theorem. This part requires a good control on the covariant derivatives of the curvature tensor and  the distance functions. This is provided by the appendices. 

\textbf{Acknowledgements.} I wish to thank Gilles Carron for drawing my attention to the geometry
of asymptotically flat manifolds and for many discussions. This work benefited from the French ANR grant GeomEinstein. 

\newpage

\tableofcontents

\section{Examples.}

\subsection{Flat plane bundles over the circle.}

To have a clear picture in mind, it is useful to understand flat manifolds obtained as quotients of the Euclidean space $\Rl^3$ by the action of a screw operation $\rho$. Let us suppose this rigid motion is the composition of a rotation of angle $\theta$ and of a translation of a distance $1$ along the rotation axis. The quotient manifold is always diffeomorphic to $\Rl^2 \times \Sph^1$, but its Riemannian structure depends on $\theta$ : one obtains a flat plane bundle over the circle whose holonomy is the rotation of angle $\theta$. These very simple examples conceals interesting features, which shed light on the link between injectivity radius, volume growth and holonomy. In this paragraph, we stick to dimension $3$ for the sake of simplicity, but what we will observe remains relevant in higher dimension.

When the holonomy is trivial, i.e. $\theta=0$, the Riemannian manifold is nothing but the standard $\Rl^2 \times \Sph^1$. The volume growth is uniformly comparable to that of the Euclidean $\Rl^2$:
$$
\Exists{A,\,B}{\Rl_+^*} \Forall{x}{M} \forall \, t \geq 1, \, A t^2 \leq \vol B(x,t) \leq B t^2.
$$
The injectivity radius is $1/2$ at each point, because of the lift of the base circle, which is even a closed geodesic; the iterates of these loops yield closed geodesics whose lengths describe all the natural integers, at each point.  

Now, consider an angle $\theta = 2\pi \cdot p/q$, for some coprime numbers $p$, $q$. A covering of order $q$ brings us back to the trivial case. The volume growth is thus uniformly comparable to that of $\Rl^2$. What about the injectivity radius ? Because of the cylindric symmetry, it depends only on the distance to the "soul", that is the image of the screw axis: let us denote by $\inj(t)$ the injectivity radius at distance $t$ from the soul. This defines a continuous function admitting uniform upper and lower bounds, but non constant in general. The soul is always a closed geodesic, so that $\inj(0)=1/2$. But as $t$ increases, it becomes necessary to compare the lengths $l_k(t)$ of the geodesic loops obtained as images of the segments $[x,\rho^k(x)]$, with $x$ at distance $t$ from the axis. We can give a formula:
\begin{equation}\label{lglacets}
l_k(t) = \sqrt{k^2 + 4 t^2 \sin^2 (k \theta/2)}.
\end{equation}
The injectivity radius is given by $2 \inj(t)=\inf_k l_k(t)$. In a neighbourhood of $0$, $2 \inj$ equals $l_1$ ; then $2 \inj$ may coincide with $l_k$ for different indices $k$. If $k < q$ is fixed, since $\sin \frac{k \theta}{2}$ does not vanish, the function $t \mapsto l_k(t)$ grows linearly and tends to infinity. The function $l_q$ is constant at $q$ and $l_q \leq l_k$ for $k \geq q$. Thus, outside a compact set, the injectivity radius is constant at $q/2$ and it is half the length of a unique geodesic loop which is in fact a closed geodesic. Besides, the other loops are either iterates of this shortest loop, or they are much longer ($l_k(t) \asymp t$). 

\begin{figure}[htb!]
\begin{center}
\input{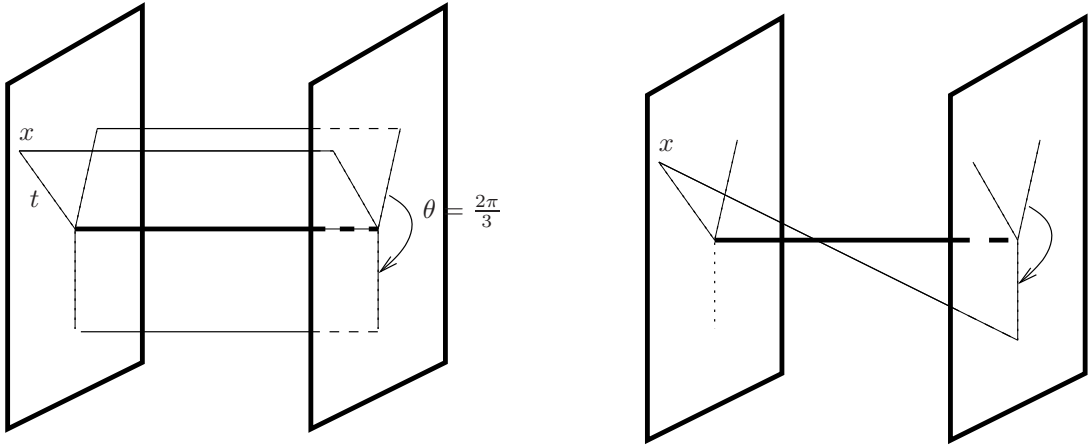}
\end{center}
\caption{The holonomy angle is $\theta=\frac{2\pi}{3}$. On the left, a geodesic loop based at $x$ with length $l_3(t)=3$. On the right, a geodesic loop based at $x$ with length $l_1(t)=\sqrt{1+9t^2}$.}
\label{fiberplat}
\end{figure}

When $\theta$ is an irrational multiple of $2\pi$, the picture is much different. In particular, the injectivity radius is never bounded. 

\begin{prop}
The injectivity radius is bounded if and only if $\theta$ is a rational multiple of $2\pi$.
\end{prop}

\proof
The "only if" part is settled, so we assume the function $t \mapsto \inj(t)$ is bounded by some number $C$. For every $t$, there is an integer $k(t)$ such that $2 \inj(t) = l_{k(t)}$. Formula (\ref{lglacets}) implies the function $t \mapsto k(t)$ is bounded by $C$. Since its values are integers, there is a sequence $(t_n)$ going to infinity and an integer $k$ such that $k(t_n) = k$ for every index $n$. Then (\ref{lglacets}) yields 
$$
\Forall{n}{\Nl} \quad  l_k(t_n)^2 = k^2 + 4 t_n^2 \sin^2 (k \theta / 2) \leq C^2.
$$   
Since $t_n$ goes to infinity, this requires $\sin^2 \frac{k \theta}{2} = 0$: there is an integer $m$ such that $k \theta/2 = m \pi$, i.e. $\theta/2\pi = m/k$.
\endproof

What about volume growth ? The volume of balls centered in some given point grows quadratically:
$$
\forall \, x, \, \exists\,  B_x, \, \forall \, t \geq 1, \, \vol B(x,t) \leq B_x t^2.
$$
In the ``rational'' case, the estimate is even uniform with respect to the center $x$ of the ball:
\begin{equation}\label{majunif}
\exists\,  B, \, \forall \, x, \, \forall \, t \geq 1, \, \vol B(x,t) \leq B t^2.
\end{equation}
In the ``irrational'' case, this strictly subeuclidean estimate is never uniform. Why ? The proposition above provides a sequence of points $x_n$ such that $r_n :=\inj(x_n)$ goes to infinity. Given a lift $\hat{x}_n$ of $x_n$ in $\Rl^3$, the ball $B(\hat{x}_n,r_n)$ is the lift of $B(x_n,r_n)$ and its volume is $\frac{4}{3} \pi r_n^3$. If we assume two points $v$ and $w$ of $B(\hat{x}_n,r_n)$ lift the same point $y$ of $B(x_n,r_n)$, there is by definition an integer number $k$ such that $\rho^k(v)=w$; since $\rho$ is an isometry  of $\Rl^3$, we get 
$$
\abs{\rho^k(\hat{x}_n) - \hat{x}_n} \leq \abs{\rho^k(\hat{x}_n) - \rho^k(v)} + \abs{\rho^k(v) - \hat{x}_n} = \abs{\hat{x}_n - v} + \abs{w - \hat{x}_n} < 2 r_n = 2 \inj(x_n),
$$ 
which contradicts the definition of $\inj(x_n)$ (the segment $[\rho^k(\hat{x}_n),\hat{x}_n]$ would go down as a too short geodesic loop at $x_n$). Therefore $B(\hat{x}_n,r_n)$ and $B(x_n,r_n)$ are isometric, hence $\vol B(x_n,r_n) = \frac{4}{3} \pi r_n^3$, which prevents an estimate like (\ref{majunif}). 

For concreteness, we wish to conclude this paragraph with quantitative estimates on the injectivity radius in the irrational case. It is indeed a Diophantine approximation problem.

\begin{prop}
Assume $\theta/(2\pi)$ is an irrational but algebraic number. Given $\alpha$ in $]0,1/2[$, there are positive numbers $T$, $C$ and $A$ depending only on $\theta$ and $\alpha$ and such that for every $t > T$, if $x$ is at distance $t$ from the soul, then 
$$
\inj(t) \geq C t^\alpha \;
\text{ and } \;
\vol B(x, C t^{\alpha}) = A t^{3\alpha}.
$$
\end{prop}

\proof
Roth theorem provides a positive constant $C = C(\theta,\alpha)$ such that for every $(k,m)$ in $\Ir^* \times I_r$:
$$
\abs{\frac{\theta}{2\pi} - \frac{m}{k}} \geq \frac{C}{k^{1/\alpha}}.
$$  
Without loss of generality, we can suppose $C \leq 1/2$. It follows that for every positive integer $k$,
$$
C k^{1-1/\alpha} \leq \frac{1}{2\pi} \abs{k \theta - 2 \pi m} \leq \frac{1}{4} \abs{e^{i k \theta} - 1},
$$  
so that we obtain for every positive number $t$ and for every integer $k$ in $[1,t^\alpha]$:
\begin{equation}\label{roth}
\abs{t e^{i k \theta} - t} \geq 4 C t^{\alpha}.
\end{equation}
Consider a ball $B:=B(x,C t^\alpha)$ centered in a point $x$ at distance $t$ from the soul. Let $\hat{x}$ be a lift of $x$ in $\Rl^3$ and let $\hat{B}$ be the ball centered in $\hat{x}$ and with radius $C t^\alpha$ in $\Rl^3$. To prove both estimates, it is sufficient to ensure this ball covers $B$ only once. Suppose $v$ and $w$ are two distinct lifts in $\hat{B}$ of a single element $y$ of $B$: $w=\rho^k v$ for some $k$. Denote by $s$ the common distance of $v$ and $w$ to the soul. The following holds:
$
\abs{w-v}^2 = k^2 + \abs{e^{i k \theta} - 1}^2 s^2.
$
Since $w$ and $v$ belong to the ball $\hat{B}$, we get 
$
k \leq \abs{w-v} \leq 2 C t^\alpha \leq t^\alpha.
$
Apply (\ref{roth}) to obtain:
$$
\abs{w-v} \geq  \abs{e^{i k \theta} - 1} s \geq 4 C t^{\alpha-1} s.
$$
With $\abs{w-v} \leq 2 C t^\alpha$, we find:
\begin{equation}\label{un}
s \leq \frac{t}{2}.
\end{equation}
Now, as $v$ belongs to $\hat{B}$, we can write
$
\abs{v - \hat{x}} \leq C t^\alpha,
$
which combines with the triangle inequality to yield:
\begin{equation}\label{deux}
s \geq t - C t^\alpha. 
\end{equation}
Comparing (\ref{un}) and (\ref{deux}), we arrive at $t \leq (2C)^{\frac{1}{1-\alpha}}$. If we suppose $t \geq T := 2 (2C)^{\frac{1}{1-\alpha}}$, it ensures $\hat{B}$ covers $B$ only once, hence the result. 
\endproof

This estimate on the injectivity radius is (quasi) optimal:

\begin{prop}
Whatever $\theta$ is, the growth of the injectivity radius is bounded from above by
$$
\forall \, t\geq 1, \, \inj(t) \leq \frac{\sqrt{1+4 \pi^2}}{2} \sqrt{t}.
$$
\end{prop}

\proof
Fix $t\geq 1$. The pigeonhole principle yields an integer $k$ in $[1,\sqrt{t}]$ such that 
$$
2 \abs{\sin k\theta/2} = \abs{e^{ik\theta}-1} \leq \frac{2\pi}{\sqrt{t}}.
$$
We deduce:
$$
2 \inj(t) \leq l_k(t) = \sqrt{k^2 + 4 t^2 \sin^2 \frac{k \theta}{2}} \leq \sqrt{1+4 \pi^2} \sqrt{t}.
$$
\endproof

When $\theta/(2\pi)$ admits good rational approximations, an almost rational behaviour can be recovered, with a slowly growing injectivity radius. For instance, if $\theta/(2\pi)$ is the Liouville number $\displaystyle{\sum_{n=1}^\infty} 10^{-n!}$, then 
$\displaystyle{\liminf_{t \To \infty}} \left(t^{-a} \inj(t) \right) = 0$
for every $a>0$.


\subsection{The Taub-NUT metric.}

The Taub-NUT metric is the basic non trivial example of ALF gravitational instanton. This Riemannian metric over $\Rl^4$ was introduced by Stephen Hawking in \cite{Haw}. To describe it, we mostly follow \cite{Leb}.

Thanks to the Hopf fibration (Chern number: $-1$), we can see $\Rl^4 \backslash \set{0} = \Rl_+^* \times \Sph^3$ as the total space of a principal circle bundle $\pi$ over $\Rl_+^* \times \Sph^2 = \Rl^3 \backslash \set{0}$. Denote by $g_{\Rl^3}$ the standard metric on $\Rl^3$ and by $r$ the distance to $0$ in $\Rl^3$. Let $V$ be the harmonic function which is defined on $\Rl^3 \backslash \set{0}$ by:
$$
V = 1 + \frac{1}{2 r}.
$$
In polar coordinates, the volume form on $\Rl^3$ reads $r^2 dr \wedge \Omega$, so that the definition
$$
d V \wedge * d V = \abs{dV}^2 r^2 dr \wedge \Omega
$$
yields $* d V = -\frac{1}{2} \Omega$. The Chern class $c_1$ of the $U(1)$-bundle $\pi$ is the cohomology class of $-\frac{\Omega}{4\pi}$. Chern-Weil theory asserts a $2$-form $\alpha$ with values in $\mathfrak{u}_1 = i \Rl \id_\Cx$ is the curvature of a connection on $\pi$ if and only if $\frac{i}{2\pi} \Tr \alpha$ represents the cohomology class $c_1$. The identity
$$
\frac{i}{2\pi} \Tr (* dV \otimes i \id) = -\frac{\Omega}{4\pi}
$$
therefore yields a connection with curvature $* dV \otimes i \id_\Cx$ on $\pi$. Let $\omega \otimes i \id_\Cx$ be the $1$-form of this connection. Lifts of objects on the base will be endowed with a hat: for instance, $\hat{V} = V \circ \pi$. 

On $\Rl^4 \backslash \set{0}$, the Taub-NUT metric is given by the formula 
$$
g= \hat{V} \hat{g}_{\Rl^3} + \hat{V}^{-1} \omega^2.
$$
Setting $\rho=\sqrt{2r}$, we obtain the behaviour of the metric near $0$: 
$$
g \cong d\rho^2 + \rho^2 \frac{\hat{g}_{\Sph^2}}{4} + \rho^2 \omega^2 
= d\rho^2 + \rho^2 \hat{g}_{\Sph^3}.
$$
It can thus be extended as a complete metric on $\Rl^4$ by adding one point, sent on the origin of $\Rl^3$ by $\pi$. The additional point can be seen as the unique fixed point for the action of $\Sph^1$.

By construction, the metric is $\Sph^1$-invariant, fiber length tends to a (nonzero) constant at infinity, while the induced metric on the base is asymptotically Euclidean(it is at distance $\OO(r^{-1})$ from the Euclidean metric). Thus there are positive constants $A$ and $B$ such that
$$
\forall\, R \geq 1, \; A R^3 \leq \vol B(x,R) \leq B R^3.
$$  

Moreover, the Taub-NUT metric is hyperk\"ahler. K\"ahler structures can be described in the following way. Let $(x,y,z)$ be the coordinates on $\Rl^3$ and choose a local gauge: $\omega = dt + \theta$ for some vertical coordinate $t$ and a $1$-form $\theta$ such that $d\theta = *dV$. An almost complex structure $J_x$ can be defined by requiring the following action on  the cotangent bundle:
$$
J_x \left(\sqrt{\hat{V}} d\hat{x} \right) = \frac{1}{\sqrt{\hat{V}}} \omega \quad \text{ et } \quad J_x \left(\sqrt{\hat{V}} d\hat{y} \right) = \sqrt{V} d\hat{z}.
$$    
It is shown in \cite{Leb} that $(g,J_x)$ is indeed a K\"ahler structure ($J_x$ being parallel outside $0$, it smoothly extends  on the whole $\Rl^4$). A circular permutation of the roles of $x$, $y$ and $z$ yields three K\"ahler structure $(g,J_x)$, $(g,J_y)$, $(g,J_z)$ satisfying $J_x J_y = J_z$, hence the hyperk\"ahler structure. \cite{Leb} even shows that these complex structures are biholomorphic to that of $\Cx^2$. Note the Taub-NUT metric can also be obtained as a hyperk\"ahler quotient (\cite{Bes}). 

In an orthonormal and left invariant trivialization $(\sigma_1, \sigma_2, \sigma_3)$ of $T^* \Sph^3$, the metric can be written
$$
g = V(r) dr^2 + 4 r^2 V(r) ( \sigma_1^2 + \sigma_2^2) + \frac{1}{V(r)} \sigma_3^2
$$
(we forget the hats). Let $H$ be the solution of 
$$
H' = \frac{1}{\sqrt{V(H)}}
$$
with $H(0)=0$ and set $r=H(t)$. Note that $H' \sim 1$ and $H \sim t$ at infinity. The equation above becomes:
$$
g = dt^2 + \left( \frac{2 H}{H'} \right)^2 ( \sigma_1^2 + \sigma_2^2) + H'^2 \sigma_3^2.
$$
Using \cite{Unn}, it is possible to compute the curvature of such a metric. It decays at a cubic rate:
$$
\abs{\Rm} = \OO(r^{-3}).
$$

This ansatz produces a whole family of examples: the "multi-Taub-NUT" metrics or $A_{N-1}$ ALF instantons. These are obtained as total spaces of a circle bundle $\pi$ over $\Rl^3$ minus some points $p_1$, $\dots$, $p_N$, endowed with the metric 
$$
\hat{V} \hat{g}_{\Rl^3} + \hat{V}^{-1} \omega^2,
$$
where $V$ is the function defined on $\Rl^3 \backslash \set{p_1,\cdots,p_N}$ by
$$
V(x) = 1 + \sum_{i=1}^N \frac{1}{2 \abs{x-p_i}}
$$
and where $\omega$ is the form of a connexion with curvature $*dV \otimes i \id_\Cx$. As above, a completion by $N$ points is possible. The circle bundle restricts on large sphere as a circle bundle of Chern number $-N$. The metric is again hyperk\"ahler and has cubic curvature decay. The underlying complex manifold is $\Cx^2/\Ir_N$. The geometry at infinity is that of the Taub-NUT metric, modulo an action of $\Ir_N$, which is the fundamental group of the end.

Other examples are built in \cite{ChH}: the geometry at infinity of these $D_k$ ALF gravitational instantons is essentially that of a quotient of a multi-Taub-NUT metric by the action of a reflection on the base.


\section{Injectivity radius and volume growth.}

\subsection{An upper bound on the injectivity radius.}

\begin{prop}[Upper bound on the injectivity radius]\label{majinj}
There is a universal constant $C(n)$ such that on any complete Riemannian manifold $(M^n,g)$ satisfying \begin{equation}\label{volume}
\inf_{t>0} \limsup_{x \To \infty} \frac{\vol B(x,t)}{t^n} < C(n),
\end{equation}
the injectivity radius is bounded from above, outside a compact set.
\end{prop}

The assumption (\ref{volume}) means there is a positive number $T$ and a compact subset $K$ of $M$ such that:
\begin{equation}\label{volumebis}
\Forall{x}{M \backslash K} \vol B(x,T) < C(n) T^n.
\end{equation}
We think of a situation where there is a function $\omega$ going to zero at infinity and such that for any point $x$, $\vol B(x,t) \leq \omega(t) t^n$. The point is we require a \emph{uniform} strictly subeuclidean volume growth. Even in the flat case, we have seen that a uniform estimate moderates the geometry much more than a centered strictly subeuclidean volume growth.

\proof
The constant $C(n)$ is given by Croke inequality \cite{Cro}:
\begin{equation}\label{croke}
\forall \, t \leq \inj(x), \, \Forall{x}{M} \vol B(x,t) \geq C(n) t^n.
\end{equation}
Let $x$ be a point outside the compact $K$ given by (\ref{volumebis}). If $\inj(x)$ is greater than the number $T$ in (\ref{volumebis}), (\ref{croke}) yields:
$$
C(n) T^n \leq \vol B(x,T) < C(n) T^n,
$$
which is absurd. The injectivity radius at $x$ is thus bounded from above by $T$.  
\endproof

Cheeger-Fukaya-Gromov theory applies naturally in this setting: it describes the geometry of Riemannian manifolds with small curvature and injectivity radius bounded from above \cite{CG}. Let us quote the

\begin{cor}
Let $(M^n,g)$ be a complete Riemannian manifold whose curvature goes to zero at infinity and satisfying (\ref{volume}). Outside a compact set, $M$ carries a $F$-structure of positive rank whose orbits have bounded diameter. 
\end{cor}

It means we already know there is some kind of structure at infinity on these manifolds. Our aim is to make it more precise, under additional assumptions.


\subsection{The fundamental pseudo-group.}

The notion of "fundamental pseudo-group" was introduced by M. Gromov in the outstanding \cite{GLP}. It is very natural tool in the study of manifolds with small curvature and bounded injectivity radius. Let us give some details. 

Let $M$ be a complete Riemannian manifold and let $x$ be a point in $M$. We assume the curvature is bounded by $\Lambda^2$ ($\Lambda \geq 0$) on the ball $B(x,2\rho)$, with $\Lambda \rho < \pi/4$. In particular, the exponential map in $x$ is a local diffeomorphism on the ball $\hat{B}(0,2\rho)$ centered in $0$ and of radius $2\rho$ in $T_x M$. The metric $g$ on $B(x,2\rho)$ thus lifts as a metric $\hat{g} := \exp_x^* g$ on $\hat{B}(0,2\rho)$. We will denote by $\Exp$ the exponential map corresponding to $\hat{g}$. 

An important fact is proved in \cite{GLP}: any two points in $\hat{B}(0,2\rho)$ are connected by a unique geodesic which is therefore minimizing; moreover, balls are strictly convex in this domain.  

When the injectivity radius at $x$ is greater than $2\rho$, the Riemannian manifolds $(B(x,\rho),g)$ and $(\hat{B}(0,\rho),\hat{g})$ are isometric. But if it is small, there are short geodesic loops based at $x$ and $x$ admits differents lifts in $\hat{B}(0,\rho)$. The fundamental pseudo-group $\Gamma(x,\rho)$ in $x$ and at scale $\rho$ measures the injectivity defect of the exponential map over $\hat{B}(0,\rho)$ \cite{GLP} : $\Gamma(x,\rho)$ is the pseudo-group consisting of all the continuous maps $\tau$ from $\hat{B}(0,\rho)$ to $T_xM$ which satisfy 
$$
\exp_x \circ \tau = \exp_x
$$ 
and send $0$ in $\hat{B}(0,\rho)$. Since $\exp_x$ is a local isometry, these maps send geodesics onto geodesics hence preserve distances :
they are isometries onto their image. In particular, they are smooth sections of the exponential map at $x$. 

\begin{figure}[htb!]
\begin{center}
\input{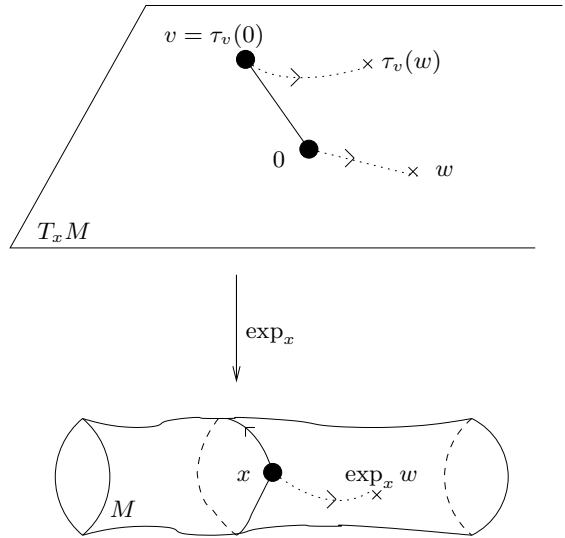}
\end{center}
\caption{$\tau_v(w)$ is obtained in the following way. Push the segment $[0,w]$ from $T_xM$ to $M$ thanks to $\exp_x$ anf lift the resulting
geodesic from $v$ to obtain a new geodesic in $T_xM$ whose tip is $\tau_v(w)$.}
\label{releve}
\end{figure}

Given a lift $v$ of $x$ in $\hat{B}(0,\rho)$ (i.e. $\exp_x(v) = p$), consider the map 
$$
\tau_v := \Exp_v \circ \left( T_v \exp_x \right)^{-1}.
$$
$\left( T_v \exp_x \right)^{-1}$ maps a point $w$ in $\hat{B}(0,\rho)$ on the initial speed vector of the geodesic lifting $t \mapsto \exp_x
tw$ from $v$. Thus $\tau_v(w)$ is nothing but the tip of this geodesic. In particular, $\tau_v(w)$ lifts $\exp_x w$, i.e. 
$\exp_x(\tau_v(w)) = \exp_x w$, 
and $\tau_v(0)$ belongs to $\hat{B}(0,\rho)$. So $\tau_v$ is an element of $\Gamma(x,\rho)$ (cf. figure \ref{releve}).   

Reciprocally, if $\tau$ is an element of $\Gamma(x,\rho)$ mapping $0$ to $v \in \hat{B}(0,\rho)$, then one can see that $\tau = \tau_v$.
Indeed, given $w$ in $\hat{B}(0,\rho)$, $t \mapsto \tau(tw)$ is the lift of $t \mapsto \exp_x tw$ from $v$, so that the argument above 
yields $\tau(w) = \tau_v(w)$.

There is therefore a one-to-one correspondance between elements of $\Gamma(x,\rho)$ and oriented geodesic loops based at $x$ with length
bounded by $\rho$. Since $\exp_x$ is a local diffeomorphism, $\Gamma(x,\rho)$ is in particular finite. Thus, given $x$, 
$(\Gamma(x,\rho))_{0<\rho<\pi/(4\Lambda)}$ is a nondecreasing finite pseudo-group family.

\begin{ex}
Consider a flat plane bundle over $\Sph^1$, with rational holonomy $\rho$ (cf. section 1): the screw angle $\theta$ is $2\pi$ times $p/q$, with
coprime $p$ and $q$. For large $\rho$ and $x$ farther than $\rho/\sin(\pi/q)$ from the soul (when $q=1$, there is no condition), the 
fundamental pseudo-group $\Gamma(x,\rho)$ is generated by the unique geodesic loop with length $q$. It therefore consists of translations
only. In particular, it does not contain $\rho$, except in the trivial case $\rho=\id$. In general, many geodesic loops are forgotten, for
they are too long.
\end{ex}

Every nontrivial element of $\Gamma(x,\rho)$ acts without fixed points. Let us show it. We assume a point $w$ is fixed by some $\tau_v$ in 
$\Gamma(x,\rho)$. First, $w$ is the tip of the geodesic $\gamma_1\, : \, t \mapsto tw$, defined on $[0,1]$ and lifting $t \mapsto \exp_x tw$ 
from $0$. Second, $w=\tau_v(w)$ is also the tip of the geodesic $\gamma_2 \, : \, t \mapsto \tau_v(tw)$ lifting the same geodesic 
$t \mapsto \exp_x tw$, but from $v$. Differentiating at $t=1$ the equality
$$
\exp_x \circ \gamma_1(t) = \exp_x \circ \gamma_2(t),
$$
we obtain $T_w \exp_x (\gamma'_1(1)) = T_w \exp_x (\gamma'_2(1))$ and thus $\gamma'_1(1)=\gamma'_2(1)$. Geodesics $\gamma_1$ and 
$\gamma_2$ must then coincide, which implies $0=\gamma_1(0)=\gamma_2(0)=v$, hence $\tau_v = \id$.

In the pseudo-group $\Gamma(x,\rho)$, every element has a well-defined inverse. To see this, given $\tau=\tau_v$ in $\Gamma(x,\rho)$, 
consider the geodesic loop $\sigma \, : \, t \mapsto \exp_x tv$ and define $\tilde{v} := - \sigma'(1)$. The maps $\tau_{\tilde{v}}$ and
$\tau_v$ are well defined on the ball $\hat{B}(0,2\rho)$, so $\tau_{\tilde{v}} \circ \tau_v$ is well defined on $\hat{B}(0,\rho)$. It is a section of $\exp_x$ and fixes $0$: $\tau_{\tilde{v}} \circ \tau_v (0) = \tau_{\tilde{v}}(v)$ 
is the tip of the geodesic lifting $\sigma$ from $\tilde{v}$; by construction, it is $0$. So $\tau_{\tilde{v}} \circ \tau_v$ is the identity. 
In the same way, one checks that $\tau_v \circ \tau_{\tilde{v}}$ is the identity. In this sense, $\tau_{\tilde{v}}$ is the inverse of 
$\tau_v$.   

Given a geodesic loop $\sigma$ with length bounded by $\rho$, let us call ``sub-pseudo-group generated by $\sigma$ in $\Gamma(x,\rho)$'' 
the pseudo-group $\Gamma_\sigma(x,\rho)$ which we describe now : it contains an element $\tau_v$ of $\Gamma(x,\rho)$ if and only if
$v$ is the tip of a piecewise geodesic segment staying in $\hat{B}(0,\rho)$ and obtained by lifting several times $\sigma$ from $0$. 
If $\tau$ is an element of $\Gamma(x,\rho)$ which corresponds to a loop $\sigma$, we will also write $\Gamma_\tau(x,\rho)$ for the
sub-pseudogroup generated by $\tau$ in $\Gamma(x,\rho)$. If $k$ is the largest integer such that $\tau^i(0)$ belongs to the ball 
$\hat{B}(0,\rho)$ for every natural number $i \leq k$, then:
$$
\Gamma_\tau(x,\rho) = \Gamma_\sigma(x,\rho) = \set{\tau^i / -k \leq i \leq k}.
$$ 
 
\begin{figure}[htb!]
\begin{center}
\input{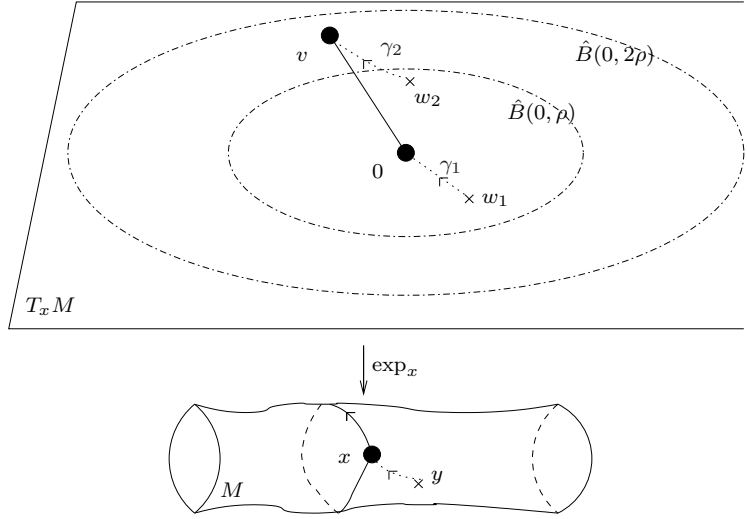}
\end{center}
\caption{$\tau_v(w_1)=w_2$.}
\label{action}
\end{figure}

If $2 \rho \leq \rho' < \frac{\pi}{4\Lambda}$, then the orbit space of the points of the ball $\hat{B}(0,\rho)$  under the action of 
$\Gamma(x,\rho')$, $\hat{B}(0,\rho)/\Gamma(x,\rho')$, is isometric to $B(x,\rho)$, through the factorization of $\exp_x$. The only thing to
check is the injectivity. Given two lifts $w_1, w_2 \in \hat{B}(0,\rho)$ of the same point $y \in B(x,\rho)$, let us prove they are 
in the same orbit for $\Gamma(x,\rho')$. Consider the unique geodesic $\gamma_1$ from $w_1$ to $0$, push it by $\exp_x$ and lift the
resulting geodesic from $w_2$ to obtain a geodesic $\gamma_2$, from $w_2$ to some point $v$ (cf. figure \ref{action}). Then $v$ 
is a lift of $x$ in $\hat{B}(0,\rho')$ (by triangle inequality) and $\tau_v$ maps $w_1$ to $w_2$, hence the result. 

We will need to estimate the number $N_x(y,\rho)$ of lifts of a given point $y$ in the ball $\hat{B}(0,\rho)$ of $T_xM$. Lifting 
one shortest geodesic loop from $0=:v_0$, we arrive at some point $v_1$. Lifting the same loop from $v_1$, we arrive at a new point $v_2$, 
etc. This construction yields a sequence of lifts $v_k$ of $x$ which eventually goes out of $\hat{B}(0,\rho)$: otherwise, since 
there cannot exist an accumulation point, the sequence would be periodic; $\tau_{v_1}$ would then fix the centre of the unique ball with 
minimal radius which contains all the points $v_k$, which is not possible, since $\tau_{v_1}$ is nontrivial hence has no fixed point 
(the uniqueness of the ball stems from the strict convexity of the balls, cf. \cite{G1}). Of course, one can do the same thing with the reverse
orientation of the same loop. Since the distance between two points $v_k$ is at least $2\inj(x)$, this yields at least $\rho/\inj(x)$ lifts
of $x$ in $\hat{B}(0,\rho)$: 
$$
\abs{\Gamma(x,\rho)} = N_x(x,\rho) \geq \rho/\inj(x).
$$
Lifting one shortest geodesic between $x$ and some point $y$ from the lifts of $x$ and estimating the distance between the tip and $0$ 
with the triangle inequality (cf. figure \ref{nbreleves}), we get: 
\begin{equation}\label{nbreleves}
N_x(y,\rho) \geq N_x(x,\rho-d(x,y)) = \abs{\Gamma(x,\rho-d(x,y))} \geq \frac{\rho-d(x,y)}{\inj(x)}.
\end{equation}
For $d(x,y)\leq \rho/2$, this yields:
\begin {equation}\label{compvol}
\frac{\rho}{2 \inj(x)} \vol B(x,\rho/2) \leq \abs{\Gamma(x,\rho/2)} \vol B(x,\rho/2) \leq \vol \hat{B}(0,\rho). 
\end {equation}

\begin{figure}[htb!]
\begin{center}
\input{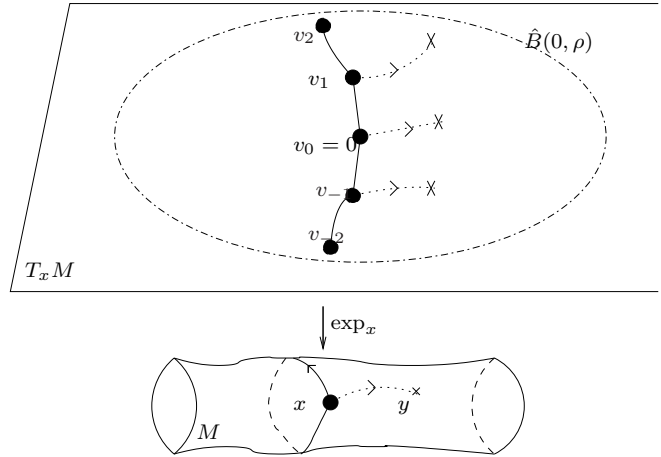}
\end{center}
\caption{Take a minimal geodesic between $x$ and $y$ and lift it from every point in the fiber of $x$ to obtain points in the fiber of $y$.}
\label{nbrelevesfig}
\end{figure}

For $\rho \leq \rho' < \frac{\pi}{4\Lambda}$, the set
$$
\mathcal{F}(x,\rho,\rho') := \set{w \in \hat{B}(0,\rho) \Big{/} \, \Forall{\gamma}{\Gamma(x,\rho')} d(0,\gamma(w)) \geq d(0,w)}
$$
is a fundamental domain for the action of $\Gamma(x,\rho')$ on the ball $\hat{B}(0,\rho)$. Finiteness ensures each orbit intersects 
$\mathcal{F}$. Furthermore, if $\tau$ belongs to $\Gamma(x,\rho')$, the set $\mathcal{F}(x,\rho,\rho') \cap \tau(\mathcal{F}(x,\rho,\rho'))$
consists of points whose distances to $0$ and $\tau(0)$ are equal, hence has zero measure: by finiteness again, up to a set with zero measure, 
$\mathcal{F}(x,\rho,\rho')$ contains a unique element of each orbit. For the same reason, if $\tau$ belongs to $\Gamma(x,\rho')$, the set
$$
\mathcal{F_\tau}(x,\rho,\rho') := \set{w \in \hat{B}(0,\rho) \Big{/} \, \Forall{\gamma}{\Gamma_\tau(x,\rho')} d(0,\gamma(w)) \geq d(0,w)}
$$
is a fundamental domain for the action of the sub-pseudo-group $\Gamma_\tau(x,\rho')$.
From our discussion follows an important fact: if $2 \rho \leq \rho' < \frac{\pi}{4\Lambda}$, then
$$
\vol \mathcal{F}(x,\rho,\rho') = \vol B(x,\rho).
$$

We will need to control the shape of these fundamental domains.

\begin{lem}\label{trigo}
Fix $\rho \leq \rho' < \frac{\pi}{4\Lambda}$ and consider a nontrivial element $\tau$ in $\Gamma(x,\rho')$. Denote by 
$\mathcal{I}_\tau(x,\rho)$ the set of points $w$ in $\hat{B}(0,\rho)$ such that 
$$
g_x\left( w,\tau(0) \right) \leq \frac{\abs{\tau(0)}^2}{2} + \frac{\Lambda^2 \rho^2 \abs{\tau(0)}^2}{2} 
$$
and
$$
g_x\left( w,\tau^{-1}(0) \right) \leq \frac{\abs{\tau(0)}^2}{2} + \frac{\Lambda^2 \rho^2 \abs{\tau(0)}^2}{2}.
$$ 
Then $\mathcal{F_\tau}(x,\rho,\rho')$ is a subset of $\mathcal{I}_\tau(x,\rho)$.
\end{lem}

\begin{figure}[htb!]
\begin{center}
\input{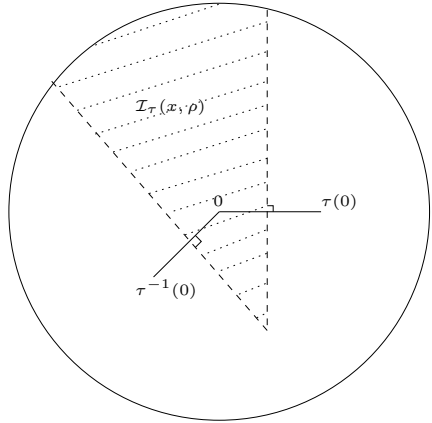}
\end{center}
\caption{The domain $\mathcal{I}_\tau(x,\rho)$.}
\label{secteur}
\end{figure}

A picture is given by figure \ref{secteur} (it represents the plane containing $0$, $\tau(0)$ and $\tau^{-1}(0)$). 

\vskip 0.5cm

\begin{figure}[htb!]
\begin{center}
\input{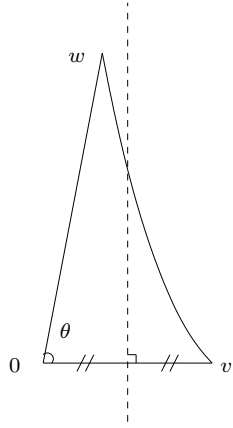}
\end{center}
\caption{$d(0,w) \leq d(v,w)$ implies $w$ is on the left of the dotted line, up to an error term.}
\label{figtrigo}
\end{figure}

\proof
Consider a point $w$ in $\mathcal{F_\tau}(x,\rho,\rho')$, set $v=\tau(0)$ and denote by $\theta \in [0,\pi]$ the angle between $v$ and $w$. 
We first assume $g_x(w,v) > 0$, that is $\theta < \pi/2$. Since any two points in $\hat{B}(0,\rho)$ are connected by a unique geodesic which
is therefore minimizing, we can apply Toponogov theorem to all triangles. In particular, in the triangle $0vw$, we find
$$
\cosh (\Lambda d(v,w)) \leq \cosh (\Lambda \abs{w}) \cosh (\Lambda \abs{v}) - \sinh (\Lambda \abs{w}) \sinh (\Lambda \abs{v}) \cos \theta.
$$
Observing $\abs{w} = d(0, w) \leq d(0,\tau^{-1} (w)) = d(v,w)$ (figure \ref{figtrigo} shows what we expect), we get 
$$
\cosh (\Lambda \abs{w}) \leq \cosh (\Lambda d(v,w)) \leq \cosh (\Lambda \abs{w}) \cosh (\Lambda \abs{v}) 
- \sinh (\Lambda \abs{w}) \sinh (\Lambda \abs{v}) \cos \theta,
$$
hence
$$
\tanh(\Lambda \abs{w}) \cos \theta \leq \frac{\cosh (\Lambda \abs{v}) - 1}{\sinh (\Lambda \abs{v})}.
$$
With $g_x(v,w) = \abs{v} \abs{w} \cos \theta$, it follows that
\begin{equation}\label{est1}
\frac{g_x(v,w)}{\abs{v}^2} \leq 
\frac{\Lambda \abs{w}}{\tanh \Lambda \abs{w}} \frac{\cosh (\Lambda \abs{v}) - 1}{\Lambda \abs{v} \sinh (\Lambda \abs{v})}.
\end{equation}
Taylor formulas yield
$$
\frac{\Lambda \abs{w}}{\tanh \Lambda \abs{w}} \leq 1 + \frac{\Lambda^2 \rho^2}{2}.
$$
and
$$
\frac{\cosh \Lambda \abs{v} - 1}{\Lambda \abs{v} \sinh \Lambda \abs{v}} \leq \frac{1}{2} + \frac{\Lambda^2 \rho^2}{6}.
$$
Combining this with $\Lambda \rho \leq \frac{\pi}{4} < 1$, one can see that (\ref{est1}) implies 
$$
\frac{g_x(v,w)}{\abs{v}^2} \leq \frac{1}{2} + \frac{\Lambda^2 \rho^2}{2}.
$$
Assuming $g_x(w,\tau^{-1}(0)) > 0$, we can work in the same way (with $v=\tau^{-1}(0)$) so as to complete the proof. 
\endproof

To understand the action of the elements in the fundamental pseudo-group, the following lemma is useful: it approximates 
them by affine transformations. 

\begin{lem}\label{holetpgf}
Consider a complete Riemannian manifold $(M,g)$ and a point $x$ in $M$ such that the curvature is bounded by $\Lambda^2$, $\Lambda \geq 0$, 
on the ball $B(x,\rho)$, $\rho>0$, with $\Lambda \rho < \pi/4$. Let $v$ be a lift of $x$ in $\hat{B}(0,\rho) \subset T_xM$. Define 
\begin{itemize}
\item the translation $t_v$ with vector $v$ in the affine space $T_xM$,
\item the parallel transport $p_v$ along $t \mapsto \exp_x tv$, from $t=0$ to $t=1$. 
\item the map $\tau_v = \Exp_v \circ \left( T_v\exp_x \right)^{-1}$,
\end{itemize}
where $\Exp$ denotes the exponential map of $(T_xM,\exp_x^*g)$. Then for every point $w$ in $\hat{B}(0,\rho-\abs{v})$,
$$
d(\tau_v(w), t_v \circ p_v^{-1}(w)) \leq \Lambda^2 \abs{v} \abs{w} (\abs{v} + \abs{w}). 
$$
\end{lem}

\proof
Proposition $6.6$ of \cite{BK} yields the following comparaison result: if $V$ is defined by $\Exp_0 V = v$ and if $W$ belongs to $T_0 T_x M$, 
then
\begin{equation}\label{bk}
d(\Exp_v \circ \hat{p}_v (W),\Exp_0(V+W)) \leq \frac{1}{3} \Lambda \abs{V} \abs{W} \sinh (\Lambda (\abs{V} + \abs{W})) \sin \angle(V,W),
\end{equation}
where $\hat{p}_v$ is the parallel transport along $t \mapsto \Exp_0 tV$, from $t=0$ to $t=1$. Set $w=\Exp_0 W$. We stress the fact that 
$Exp_0 = T_0 \exp_x$ is nothing but the canonical identification between the tangent space $T_0T_xM$ to the vector space $T_xM$ and the 
vector space itself, $T_xM$. In particular, $\Exp_0(V+W) = v+w = t_v(w)$. Since $\exp_x$ is a local isometry, we have 
$$
\hat{p}_v = \left( T_v\exp_x \right)^{-1} \circ p_v \circ T_0\exp_x, 
$$ 
so that $\Exp_v \circ \hat{p}_v (W) = \tau_v \circ p_v (w)$. With 
$$
\sinh (\Lambda (\abs{V} + \abs{W})) \leq \Lambda (\abs{V} + \abs{W}) \cosh(1) \leq 3 \Lambda (\abs{V} + \abs{W}),
$$
it follows from (\ref{bk}) that: 
$$
d(\tau_v \circ p_v (w),t_v(w)) \leq \Lambda^2 \abs{v} \abs{w} (\abs{v} + \abs{w}).
$$
Changing $w$ into $p_v^{-1} (w)$, we obtain the result.
\endproof

\subsection{Fundamental pseudo-group and volume.}

\subsubsection{Back to the injectivity radius.}

Our discussion of the fundamental pseudo-group enables us to recover a result of \cite{CGT}.

\begin{prop}[Lower bound for the injectivity radius.]\label{mininj}
Let $(M^n,g)$ be a complete Riemannian manifold. Assume the existence of $\Lambda \geq 0$ and $V>0$ such that for every point $x$ in $M$, 
$$
\abs{\Rm_x} \leq \Lambda^2 \quad \text{and} \quad \vol B(x,1) \geq V. 
$$
Then the injectivity radius admits a positive lower bound $I = I(n,\Lambda,V)$.
\end{prop}

\proof
Set $\rho = \min(1,\frac{\pi}{8\Lambda})$ and assume there is a point $x$ in $M$ and a geodesic loop based at $x$ with length bounded by 
$\rho$. Apply (\ref{compvol}) to find
$$
\frac{\rho}{2\inj(x)} \vol B(x,\rho/2) \leq \vol \hat{B}(0,\rho). 
$$
Bishop theorem estimates the right-hand side by $\omega_n \cosh(\Lambda \rho)^{n-1} \rho^n$, where $\omega_n$ is the volume of the unit ball
in $\Rl^n$. We thus obtain $\inj(x) \geq C(n,\Lambda) \vol B(x,\rho/2)$ for some $C(n,\Lambda)>0$. Since Bishop theorem also yields
a constant $C'(n,\Lambda)>0$ such that 
$$
\vol B(x,1) \leq  C'(n,\Lambda)^{-1} \vol B(x,\rho/2),
$$ 
we are left with $\inj(x) \geq C(n,\Lambda) C'(n,\Lambda) \vol B(x,1) \geq C(n,\Lambda) C'(n,\Lambda) V$.
\endproof

Combining propositions \ref{majinj} and \ref{mininj}, we obtain

\begin{cor}[Injectivity radius pinching.]\label{inj}
Let $(M^n,g)$ be a complete Riemannian manifold with bounded curvature. Suppose:
$$
\Forall{x}{M} V \leq \vol B(x,t) \leq \omega(t) t^n
$$
for some positive number $V$ and some function $\omega$ going to zero at infinity. Then there are positive numbers 
$I_1$, $I_2$ such that for any point $x$ in $M$: 
$$
I_1 \leq \inj(x) \leq I_2.
$$
\end{cor}

\subsubsection{Self-improvement of volume estimates.}

Here and in the sequel, we will always distinguish a point $o$ in our complete non-compact Riemannian
mainfolds. The distance function to $o$ will always be denoted by $r_o$ or $r$. We will always assume our manifolds 
are smooth and connected. 

\begin{prop}\label{majvol}
Let $(M^n,g)$ be a complete non-compact Riemannian manifold with faster than quadratic curvature decay, i.e.
$$
\abs{\Rm} = \OO(r^{-2-\epsilon})
$$
for some $\epsilon>0$. If there exists a function $\omega$ which goes to zero at infinity and satisfies
$$
\Forall{x}{M} \forall t \geq 1, \, \vol B(x,t) \leq \omega(t) t^n,
$$
then there is in fact a number $B$ such that 
$$
\Forall{x}{M} \forall t \geq 1, \, \vol B(x,t) \leq B t^{n-1}.
$$
\end{prop}

Indeed, under faster than quadratic curvature decay assumption, if $C(n)$ denotes the constant in \ref{majinj}, the property
$$
\inf_{t>0} \limsup_{x \To \infty} \frac{\vol B(x,t)}{t^n} < C(n),
$$
automatically implies 
$$
\limsup_{t \To \infty} \sup_{x \in M} \frac{\vol B(x,t)}{t^{n-1}} < \infty.
$$

\proof
Proposition \ref{majinj} yields an upper bound $I_2$ on the injectivity radius. Our assumption on the curvature 
implies that, given a point $x$ in $M \backslash B(o,R_0)$, with large enough $R_0$, one can apply (\ref{compvol}) with 
$2 I_2 \leq \rho = 2 t \leq r(x)/2$:
$$
\frac{t}{\inj(x)} \vol B(x,t) \leq \vol \hat{B}(0,2t). 
$$
Thanks to the curvature decay, if $R_0$ is large enough, Bishop theorem bounds the right-hand side by $\omega_n \cosh(1)^{n-1} (2t)^n$; 
with proposition \ref{majinj}, it follows that for $I_2 \leq t \leq r(x)/2$:
$$
\vol B(x,t) \leq \omega_n \cosh(1)^{n-1} 2^n I_2 t^{n-1}.
$$
We have found some number $B_1$ such that for every $x$ outside the ball $B(o,R_0)$ and for every $t$ in $[I_2,r(x)/2]$,
\begin{equation}\label{eq1}
\vol B(x,t) \leq B_1 t^{n-1}.
\end{equation}
From lemma 3.6 in \cite{LT}, which refers to the construction in the fourth paragraph of \cite{A2}, we can find a number $N$ 
such that for any natural number $k$, setting $R_k=R_0 2^k$, the annulus $A_k:=B(o,2 R_k) \backslash B(o, R_k)$
is covered by a family of balls $(B(x_{k,i}, R_k/2))_{1 \leq i \leq N}$ centered in $A_k$. Since the volume of the balls 
$B(x_{k,i},R_k/2)$ is bounded by $B_1 (R_k/2)^{n-1}$, we deduce the existence of a constant $B_2$ such that for every $t \geq I_2$,  
$$
\vol B(o,t) \leq B_2 \sum_{k=0}^{\lceil \log_2 (t/R_0) \rceil} (2^k)^{n-1},
$$
and thus, for some new constant $B_3$, we have    
\begin{equation}\label{eq2}
\forall \, t \geq I_2, \, \vol B(o,t) \leq B_3 \, t^{n-1}.
\end{equation}
Now, for every $x$ in $M \backslash B(o,R_0)$ and every $t\geq r(x)/4$, we can write 
$$
\vol B(x,t) \leq \vol B(o,t+r(x)) \vol B(o,5t) \leq 5^{n-1} B_3 t^{n-1}.
$$
And when $x$ belongs to $B(o,R_0)$, for $t \geq I_2$, we observe
$$
\vol B(x,t) \leq \vol B(o,t+R_0) \leq \vol B(o,(1+R_0/2) t)  \leq B_3 (1+R_0/2)^{n-1} t^{n-1}.
$$
Therefore there is a constant $B$ such that for every $x$ in $M$ and every $t \geq I_2$, the volume 
of the ball $B(x,t)$ is bounded by $B t^{n-1}$. The result immediately follows.
\endproof

When the Ricci curvature is nonnegative, the assumption on the curvature can be relaxed. 

\begin{prop}
Let $(M^n,g)$ be a complete non-compact Riemannian manifold with nonnegative Ricci curvature 
and quadratic curvature decay, i.e.
$$
\abs{\Rm} = \OO(r^{-2}).
$$
If there exists a function $\omega$ which goes to zero at infinity and satisfies
$$
\Forall{x}{M} \forall t \geq 1, \, \vol B(x,t) \leq \omega(t) t^n,
$$
then there is in fact a number $B$ such that 
$$
\Forall{x}{M} \forall t \geq 1, \, \vol B(x,t) \leq B t^{n-1}.
$$
\end{prop}

\proof
The previous proof can easily be adapted. (\ref{eq1}) holds for $I_2 \leq t \leq \delta r(x)$, with a small $\delta>0$. The existence 
of the covering leading to (\ref{eq2}) stems from Bishop-Gromov theorem (the $x_{k,i}$ are given by a maximal $R_k/2$-net). 
\endproof

This threshold effect shows that the first collapsing situation to study is that of a ``codimension $1$'' collapse, where the volume 
of balls with radius $t$ is (uniformly) comparable to $t^{n-1}$. This explains the gap between ALE and ALF gravitational instantons, under a uniform upper bound on the volume growth: there is no gravitational instanton with intermediate volume growth, between $\vol B(x,t) \asymp t^3$ and 
$\vol B(x,t) \asymp t^4$.

\section{Collapsing at infinity.}

\subsection{Local structure at infinity.}

We turn to codimension 1 collapsing at infinity. It turns out that the holonomy of 
short geodesic loops plays an important role. In order to obtain a nice structure, we will make a strong assumption on it. The next paragraph will explain why gravitational instantons satisfy this assumption. 

\begin{prop}[Fundamental pseudo-group structure]\label{pgf}
Let $(M^n,g)$ be a complete Riemannian manifold such that
$$
\forall \, x \in M, \, \forall \, t \geq 1, \, A t^{n-1} \leq \vol B(x,t) \leq B t^{n-1}
$$
with $0 < A \leq B$. We assume there is constant $c > 1$ such that 
$$
\abs{\Rm} \leq c^2 r^{-2}
$$
and such that if $\gamma$ is a geodesic loop based at $x$ and with length $L \leq c^{-1} r(x)$, then the holonomy $H$ of $\gamma$ satisfies
$$
\abs{H - \id} \leq \frac{c L}{r(x)}.
$$  
Then there exists a compact set $K$ in $M$ such that for every $x$ in $M \backslash K$, there is a unique geodesic loop $\sigma_x$ of minimal length $2\inj(x)$. Besides there are geometric constants $L$ and $\kappa > 0$  such that the fundamental pseudo-group $\Gamma(x,\kappa r(x))$ has at most $L r(x)$ elements, all of which are obtained by successive lifts of $\sigma$.  
\end{prop}

\begin{defn}
$\sigma_x$ is the ``fundamental loop at $x$''. 
\end{defn}

\proof
Let us work around a point $x$ far away from $o$, say with $r(x) > 100 I_2 c$. Recall (\ref{inj}) yields constants $I_1$, $I_2$ such that $0<I_1 \leq \inj \leq I_2$. The fundamental pseudo-group $\Gamma := \Gamma(x,\frac{r(x)}{4c})$ contains the sub-pseudo-group $\Gamma_\sigma :=\Gamma_\sigma(x,\frac{r(x)}{4c})$ corresponding to the loop $\sigma$ of minimal length $2\inj(x)$. Denote by $\tau=\tau_v$ one of 
the two elements of $\Gamma$ that correspond to $\sigma$: $\abs{v}=2\inj(x)$. (\ref{compvol}) implies for $\rho= \frac{r(x)}{2c}$:
$$
\abs{\Gamma} \vol B\left(x, \frac{r(x)}{4c} \right) \leq \vol \hat{B}\left(0, \frac{r(x)}{2c} \right). 
$$
Bishop theorem bounds (from above) the Riemannian volume of $\hat{B}(0, \frac{r(x)}{2c})$ by $\left( \cosh c \right)^n$ times its Euclidean volume. With the lower bound on the volume growth, we thus obtain: 
$$
\abs{\Gamma} A  \left( \frac{r(x)}{4 c} \right)^{n-1} \leq \left( \cosh c \right)^n \omega_n \left( \frac{r(x)}{2 c} \right)^n, 
$$
where $\omega_n$ denotes the volume of the unit ball in $\Rl^n$. We deduce the estimate 
$$
\abs{\Gamma} \leq L r(x)
$$ 
with 
$$
L :=\frac{2^{n-2} \omega_n \left( \cosh c \right)^n }{A c}.
$$
Now, consider an oriented geodesic loop $\gamma$, based at $x$ and with length inferior to $\frac{r(x)}{4 c}$. Name $\tau_z$ the corresponding element of $\Gamma := \Gamma(x,\frac{r(x)}{4 c})$. $H_z$ will denote the holonomy of the opposite orientation of $\gamma$. By assumption,
$$
\abs{H_z - \id} \leq \frac{c \abs{z}}{r(x)}.  
$$
The vector $z=\tau_z(0)$ is the initial speed of the geodesic $\gamma$, parameterized by $[0,1]$ in the chosen orientation. In the same way,  $\tau_z^{-1}(0)$ is the initial speed vector of $\gamma$, parameterized by $[0,1]$, but in the opposite orientation.  We deduce $-z$ is obtained as the parallel transport of $\tau_z^{-1}(0)$ along $\gamma$: $H_z (\tau_z^{-1}(0)) = -z$. From the estimate above stems:
\begin{equation}\label{plusmoins}
\abs{\tau_z^{-1}(0) + z} \leq \frac{c \abs{z}^2}{r(x)}.
\end{equation}
Given a small $\lambda$, say $\lambda = \frac{1}{100 c}$, we consider a point $w$ in the domain $\mathcal{I}_{\tau_z}(x,\lambda r(x))$ (see the definition in \ref{trigo}). It satisfies 
$$
g_x(w,\tau_z^{-1}(0)) \leq \frac{\abs{z}^2}{2} + 2 c^2  \lambda^2 \abs{z}^2.
$$
With  
$$
g_x(w,z) = -g_x(w,\tau_z^{-1}(0)) + g_x(w,\tau_z^{-1}(0)+z) 
\geq -g_x(w,\tau_z^{-1}(0)) - \abs{w} \abs{\tau_z^{-1}(0) + z},
$$
we find
$$
g_x(w,z) \geq - \frac{\abs{z}^2}{2} - 2 c^2  \lambda^2 \abs{z}^2 -  \lambda c \abs{z}^2 , 
$$
that is
$$
g_x(w,z) \geq - \frac{\abs{z}^2}{2} \left( 1 + 4 c^2 \lambda^2 + 2 \lambda c \right). 
$$
With lemma \ref{trigo}, this ensures: 
$$
\mathcal{F}_{\tau_z}\left(x,\lambda r(x), \frac{r(x)}{4}\right) \subset 
\set{w \in \hat{B}(0,\lambda r(x)) \, \Big{/} \, \abs{g_x\left( w,z \right)} \leq \frac{\abs{z}^2}{2} \left( 1 + 4 c^2 \lambda^2 + 2 \lambda c \right) }.
$$
And with $\lambda = \frac{1}{100 c}$, this leads to   
\begin{equation}\label{controledom}
\mathcal{F}_{\tau_z}\left(x,\lambda r(x), \frac{r(x)}{4}\right) \subset 
\set{w \in \hat{B}(0,\lambda r(x)) \, \Big{/} \, \abs{g_x\left( w,z \right)} \leq \frac{3\abs{z}^2}{4}}.
\end{equation}

Let $\tau'$ be an element of $\Gamma \backslash \Gamma_\sigma$ such that $v':=\tau'(0)$ has minimal norm. Suppose $\abs{v'} < \lambda r(x)$. Then, 
the minimality of $\abs{v'}$ combined with (\ref{controledom}) yields
$$
\abs{g_x\left( v',v \right)} \leq \frac{3\abs{v}^2}{4}.
$$
If $\theta \in [0,\pi]$ is the angle between $v$ and $v'$, this means:
$\abs{v'} \abs{\cos \theta} \leq 0.75 \abs{v}$. Since $\abs{v} \leq \abs{v'}$, we deduce $\abs{\cos \theta} \leq 0.75$, hence $\sin \theta \geq 0.5$. Applying (\ref{controledom}) to $\tau$ and $\tau'$, we also get
\begin{eqnarray*}
\mathcal{F}\left(x,\lambda r(x),\frac{r(x)}{4c}\right) &\subset& 
\mathcal{F}_{\tau_v}\left(x,\lambda r(x),\frac{r(x)}{4c}\right) \cap \mathcal{F}_{\tau_{v'}}\left(x,\lambda r(x),\frac{r(x)}{4c}\right)  \\
&\subset& \set{w \in \hat{B}(0,\lambda r(x)) \, \Big{/} \, \abs{g_x\left( w,v \right)} \leq \abs{v}^2, \; \abs{g_x\left( w,v' \right)} \leq \abs{v'}^2 }.
\end{eqnarray*}

\begin{figure}[htb!]
\begin{center}
\input{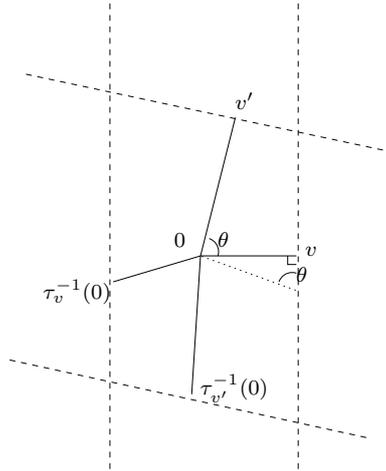}
\end{center}
\caption{The fundamental domain is inside the dotted line.}
\label{domaine}
\end{figure}

The Riemannian volume of $\mathcal{F}(x,\lambda r(x),r(x)/(4c))$ equals that of $B(x,\lambda r(x))$, so it is not smaller than $A \lambda^{n-1} r(x)^{n-1}$. The Euclidean volume of 
$$
\set{w \in \hat{B}(0,\lambda r(x)) \, \Big{/} \, \abs{g_x\left( w,v \right)} \leq \abs{v}^2 \text{ and } \abs{g_x\left( w,v' \right)} \leq \abs{v'}^2 }
$$
is not greater than $4\abs{v} \abs{v'} (2 \lambda r(x))^{n-2} / \sin \theta \leq 2^{n+2} \lambda^{n-2} I_2 \abs{v'} r(x)^{n-2}$. Comparison yields 
$$
A \lambda^{n-1} r(x)^{n-1} \leq 2^{n+2} \left(\cosh c \right)^n \lambda^{n-2} I_2 \abs{v'} r(x)^{n-2}, 
$$
that is
$$
\abs{v'} \geq \frac{\lambda A}{2^{n+2} I_2 \left(\cosh c \right)^n} r(x).
$$
Given a positive number $\kappa$ which is smaller than $\lambda$ and $\frac{\lambda A}{2^{n+2} I_2 \left(\cosh c \right)^n}$, we conclude that for any $x$ outside some compact set, the pseudo-group $\Gamma(x,\kappa r(x))$ only consists of iterates of $\tau$ (in $\Gamma(x,\frac{r(x)}{4c})$).

Suppose there are two geodesic loops with minimal length $2 \inj(x)$ at $x$. They correspond to distinct points $v$ and $v'$ in $T_xM$. Write
\begin{eqnarray*}
\mathcal{F}\left(x,\lambda r(x),\frac{r(x)}{4} \right) &\subset& 
\mathcal{F}_{\tau_v}\left(x,\lambda r(x),\frac{r(x)}{4} \right) \cap \mathcal{F}_{\tau_{v'}}\left(x,\lambda r(x),\frac{r(x)}{4} \right)  \\
&\subset& \set{w \in \hat{B}(0,\lambda r(x)) \, \Big{/} \, \abs{g_x\left( w,v \right)} \leq \abs{v}^2, \; \abs{g_x\left( w,v' \right)} \leq \abs{v}^2 }.
\end{eqnarray*}
As above, we find
$$
A \lambda^{n-1} r(x)^{n-1} \leq 2^n \left(\cosh c \right)^n \lambda^{n-2} \abs{v} \abs{v'} r(x)^{n-2}/ \sin \theta, 
$$
where $\theta \in [0,\pi]$ is the angle between the vectors $v$ and $v'$. Here, $\abs{v}=\abs{v'} \leq 2 I_2$, so 
$$
A \lambda r(x) \sin \theta \leq 2^{n+2} I_2^2 \left(\cosh c \right)^n. 
$$
The minimality of $\abs{v}$ and distance comparison yield 
$$
\abs{v} \leq d(v,v') \leq \cosh (0.02) \abs{v-v'}
$$
hence $\cos \theta \leq 0.51$. For the same reason, we find
$$
\abs{v} \leq d(\tau_v^{-1}(0),v') \leq \cosh (0.02) \abs{\tau_v^{-1}(0) - v'}.
$$
With (\ref{plusmoins}), which gives
$$
\abs{\tau_v^{-1}(0) + v} \leq 0.01 \abs{v},
$$
we deduce 
$$
\abs{v+v'} \geq \abs{\tau_v^{-1}(0) - v'} - \abs{\tau_v^{-1}(0) + v} 
\geq 0.98 \abs{v}
$$
hence $\cos \theta \geq - 0.52$, then $\abs{\cos \theta} \leq 0.52$, and 
$\sin \theta \geq 0.8$. Eventually, we obtain
$$
0.8 A \lambda r(x)  \leq 2^{n+2} I_2^2 \left(\cosh c \right)^n, 
$$
which cannot hold if $x$ is far enough from $o$. This proves the uniqueness of the shortest geodesic loop.
\endproof

Uniqueness implies smoothness:

\begin{lem}\label{reglacet}
In the setting of proposition \ref{pgf}, there are smooth local parameterizations for the family of loops $(\sigma_x)_x$. More precisely, given an orientation of $\sigma_x$, we can lift it to $T_xM$ through $\exp_x$ ; denoting the tip of the resulting segment by $v$, if $w$ is in neighborhood of $0$ in $T_xM$, then the fundamental loop at $\exp_x w$ is the image by $\exp_x$ of the unique geodesic connecting $w$ to $\tau_v(w)$.  
\end{lem}

\proof
We first prove continuity. Let $y$ be in $M$ (outside the compact set $K$) and let $(y_n)$ be a sequence converging to $y$. Let $V_n$ be a sequence of of initial unit speed vector for $\sigma_{y_n}$. Compactness ensures $V_n$ can be assumed to converge to $V$. Let $\alpha$ be the geodesic emanating from $y$ with initial speed $V$. For every index $n$, we have $\exp_{y_n} (\inj(y_n) V_n) = \sigma_{y_n}(\inj(y_n)) = y_n$. 
Continuity of the injectivity radius (\cite{GLP}) allows to take a limit: $\alpha(\inj(y)) = \exp_y \inj(y) V = y$.
Uniqueness implies $\alpha$ parameters $\sigma_y$. This yields the continuity of $(\sigma_x)_x$. Now, given $w$ in a neighborhood of $0$ in $T_xM$, consider the $e(w)$ of the lift of $\sigma_{\exp_x w}$. The map $e$ is a continuous section of $\exp_x$ and $e(0)=\tau_v(0)$ : $e=\tau_v$. The result follows.
\endproof

Now we turn to gravitational instantons: we can control their holonomy and thus apply the previous proposition. 

\subsection{Holonomy in gravitational instantons.}

\begin{lem}\label{hyperhol}
Let $(M^4,g)$ be a complete hyperk\"ahler manifold with: 
$$
\inj(x) \geq I_1 > 0
\quad \text{and} \quad 
\abs{\Rm} \leq Q  r^{-3}.
$$
Then there is some positive $c=c(I_1,Q)$ such that the holonomy $H$ of 
geodesic loops based at $x$ and with length $L \leq r(x)/c$ satisfies
$$
\abs{H - \id} \leq \frac{c}{r(x)}. 
$$
\end{lem}

\proof
Consider a point $x$ (far from $o$) and a geodesic loop based at $x$, with $L \leq r(x)/4$. Let $\tau_v \in \Gamma(x,r(x)/4)$ be a corresponding element. Thanks to (\ref{holetpgf}), we know that for every point $w$ in $T_xM$ such that $\abs{w} \leq r(x)/4$:
$$
d(\tau_v(w), t_v \circ p_v^{-1}(w)) \leq 8 Q r(x)^{-3} \abs{v} \abs{w} (\abs{v} + \abs{w})
$$
and therefore
$$
d(\tau_v(w), t_v \circ p_v^{-1}(w)) \leq Q L r(x)^{-1}.
$$
Set $H = p_v^{-1}$:  
$$
d(\tau_v(w), H w + v) \leq Q L r(x)^{-1}.
$$
Since we are working on a hyperk\"ahler $4$-manifold, the holonomy group is included in $SU(2)$, so that in some orthonormal basis
of $T_xM$, seen as complex $2$-space, $H$ reads
$$
H = \left(
\begin{array}{cc}
e^{i\theta} & 0 \\
0 & e^{-i\theta}
\end{array}
\right)
$$
with an angle $\theta$ in $]-\pi,\pi]$. Suppose $\theta$ is not zero (otherwise the statement is trivial). The equation $Hw + v =w$ admits a solution:
$$
w = \left(
\begin{array}{c}
\frac{v_1}{1-e^{i\theta}} \\
\frac{v_2}{1-e^{-i\theta}}
\end{array}
\right)
$$
where $v_1$ and $v_2$ denote the coordinates of $v$. If $\abs{w} \leq r(x)/4$, we obtain
$$
d(\tau_v(w), w) \leq Q L r(x)^{-1}.
$$
The lower bound on the injectivity radius yields
$$
d(\tau_v(w) , w) \geq 2 I_1.
$$
So we find
$$
L \geq \frac{2 I_1 r(x)}{Q}.
$$
As a consequence, if $L < \frac{2 I_1 r(x)}{Q}$, then 
$$
\abs{w} = \frac{L}{\abs{1-e^{i \theta}}}  > \frac{r(x)}{4}, 
$$
that is $\abs{H - \id} = \abs{1-e^{i \theta}} \leq 4 L r(x)^{-1}$.
\endproof

As a result, we obtain the

\begin{prop}
Let $(M^4,g)$ be a complete hyperk\"ahler manifold with
$$
\int_M \abs{\Rm}^2 r dvol < \infty
$$
and  
$$
\forall \, x \in M, \, \forall \, t \geq 1, \, A t^3 \leq \vol B(x,t) \leq B t^3
$$
with $0 < A \leq B$. Then there exists a compact set $K$ in $M$ such that for every $x$ in $M \backslash K$, there is a unique geodesic loop $\sigma_x$ of minimal length $2\inj(x)$. Besides there are geometric constants $L$ and $\kappa > 0$  such that the fundamental pseudo-group $\Gamma(x,\kappa r(x))$ has at most $L r(x)$ elements, all of which are obtained by successive lifts of $\sigma$.  
\end{prop}

\proof
Since $M$ is hyperk\"ahler, it is Ricci flat. So \cite{Min} applies (see appendix A):   
$\abs{\Rm} =\OO( r^{-3})$. So we use proposition \ref{pgf}, thanks to lemma \ref{hyperhol}.  
\endproof

\begin{rem}
From now on, we will remain in the setting of four dimensional hyperk\"ahler manifolds. It should nonetheless be noticed that the only reason for this is lemma \ref{hyperhol}.  If the conclusion of this lemma is assumed and if we suppose convenient estimates on the covariant derivatives of the curvature tensor, then we can work in any dimension (see \ref{thmgl} below). 
\end{rem}


\subsection{An estimate on the holonomy at infinity.}

To go on, we will need a better estimate of the holonomy of short loops. This is the goal of this paragraph. First, let us state an easy lemma, adapted from \cite{BK}.

\begin{lem}[Holonomy comparison]\label{compholo}
Let $\appli{\gamma}{[0,L]}{N}$ be a curve in a Riemannian manifold $N$ and let $t \mapsto \alpha_t$ be a family of loops, parameterized by $0 \leq s \leq l$ with $\alpha_t(0) = \alpha_t(l) =\gamma(t)$. We denote by $p_\gamma(t)$ the parallel transportation along $\gamma$, from $\gamma(0)$ to $\gamma(t)$. We consider a vector field $(s,t) \mapsto X(s,t)$ along the family $\alpha$ and we suppose it is parallel along each loop $\alpha_t$ ($\nabla_s X(s,t) =0$) and along $\gamma$ ($\nabla_t X(0,t) =0$). Then:
$$
\abs{p_\gamma(L)^{-1}  X(l,L) -  X(l,0)}
\leq \int_0^L \int_0^l \abs{\Rm(\partial_s \sigma_t, \partial_t \sigma_t) X(s,t)} ds dt. 
$$
\end{lem}

\begin{figure}[htb!]
\begin{center}
\input{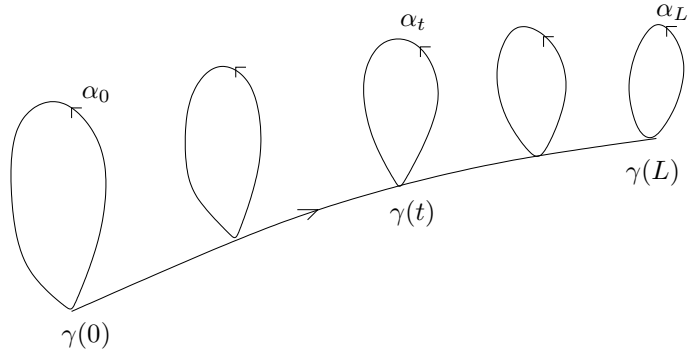}
\end{center}
\caption{A one parameter family of loops.}
\label{tube}
\end{figure}

Consider a complete hyperk\"ahler manifold $(M^4,g)$ with 
$$
\forall \, x \in M, \, \forall \, t \geq 1, \, A t^3 \leq \vol B(x,t) \leq B t^3
$$
($0 < A \leq B$) and
$$
\int_M \abs{\Rm}^2 r dvol < \infty
\quad \text{or equivalently} \quad \abs{\Rm} = \OO(r^{-3}).
$$

We choose a unit ray $\appli{\gamma}{\Rl_+}{M}$ starting from $o$ and we denote by $p_\gamma(t)$ the parallel transportation along $\gamma$, from $\gamma(0)$ to $\gamma(t)$. For large $t$, we can define the holonomy endomorphism $H_{\gamma(t)}$ of the fundamental loop $\sigma_{\gamma(t)}$: here, there is an implicit choice of orientation for the loops $\sigma_{\gamma(t)}$, which we can assume continuous. This yields an  element of $O(T_{\gamma(t)} M)$. Holonomy comparison lemma \ref{compholo} asserts that for large $t_1 \leq t_2$:
\begin{eqnarray*}
& & \abs{p_\gamma(t_2)^{-1} H_{\gamma(t_2)} p_\gamma(t_2) - p_\gamma(t_1)^{-1} H_{\gamma(t_1)} p_\gamma(t_1)} \\
&\leq& 
\int_{t_1}^\infty \int_0^1 \abs{\Rm}(c(t,s)) \abs{\partial_s c(t,s) \wedge \partial_t c(t,s))} ds dt,
\end{eqnarray*}
where, for every fixed $t$, $c(t,.)$ parameterizes $\sigma_{\gamma(t)}$ by $[0,1]$, at speed $2 \inj(\gamma(t))$.

\begin{lem}
$\abs{\partial_s c \wedge \partial_t c)}$ is uniformly bounded. 
\end{lem}

\proof
The upper bound on the injectivity radius bounds $\abs{\partial_s c}$. We need to bound the composant of $\partial_t c$ that is orthogonal to $\partial_s c$. We concentrate on a neighborhhood of some point $x$ along $\gamma$. For convenience, we change the parameterization so that $x=\gamma(0)$. We also lift the problem to $T_xM =: E$, endowed with the lifted metric $\hat{g}$. If $v=\gamma'(0)$, $\gamma$ lifts as a curve $\hat{\gamma}$ parameterized by $t \mapsto t v$. The lift $\hat{c}$ of $c$ consists of the geodesics $\hat{c}(t,.)$ connecting $t v$ to $\tau(t v)$; $\tau$ is the element of the fundamental pseudo-group corresponding to $\sigma_x$, for the chosen orientation. Observe
$$
\hat{c}(t,s) = \Exp_{t v} s X(t)
$$
where $X(t) \in T_{tv} E$ is defined by 
$$
\Exp_{t v} X(t) = \tau(t v). 
$$
The vector field $J$ defined along $\hat{c}(0,.)$ by 
$$
J(s) = \partial_t \hat{c} (0,s) = \derivzero{t} \Exp_{t v} s X(t)
$$
is a Jacobi field with initial data 
$$
J(0) = \derivzero{t} t v = v
$$  
(we identify $E$ to $T_0 E$ thanks to $\Exp_0$, in the natural way) and
\begin{eqnarray*}
J'(0) &=& \left(\nabla_s \partial_t \hat{c}  \right)\Big\vert_{(t,s)=(0,0)} \\
&=& \left(\nabla_t \partial_s \hat{c}  \right)\Big\vert_{(t,s)=(0,0)} \\
&=& \left(\nabla_t \partial_s \Exp_{t v} s X(t)  \right)\Big\vert_{(t,s)=(0,0)} \\
&=& \left(\nabla_t X \right)(0). 
\end{eqnarray*}
Suppose the curvature is bounded by $\Lambda^2$, $\Lambda >0$, in the area under consideration and apply lemma 6.3.7 of \cite{BK}: the part
$\tilde{J}$ of $J$ that is orthogonal to $\hat{c}(0,.)$ satisfies
$$
\abs{\tilde{J}(s) - p(sv) \tilde{J}(0) - s p(sv) \tilde{J}'(0)} \leq a(s)
$$
where $p(.)$ is the radial parallel transportation and where $a$ solves
$$
a'' - \Lambda^2 a = \Lambda^2 \left( \abs{\tilde{J}(0)} + \abs{\tilde{J}'(0)} \right)
$$
with $a(0)=a'(0)=0$, i.e. 
$$
a(s) = \left( \abs{\tilde{J}(0)} + \abs{\tilde{J}'(0)} \right) \left( \cosh(\Lambda s) -1 \right).
$$ 
Since here $0\leq s \leq1$ and $\Lambda << 1$, we only need a bound on $\abs{\tilde{J}(0)}$ and $\abs{\tilde{J}'(0)}$ to control $\tilde{J}$
and end the proof. Since $J(0)=v$ has unit length, we are left to bound $\tilde{J}'(0)$.

We consider the family of vectors $Y$ in $T_0 E \cong E$ that is defined by
$$
X(t)= p(tv)Y(t).
$$
With 
$$
\nabla_t X (t) = p(t v) \frac{d}{dt} p(t v)^{-1} X(t) = p(t v) Y'(t),
$$
we see that $J'(0)$ is exactly $Y'(0)$. Let $f$ be the map from $E \times T_0 E \cong E^2$ to $E$ that is defined by 
$$
f(w,W)= \Exp_w p(w) W.
$$
The equality 
$$
f(tv,Y(t))= \tau(t v).
$$
can be differentiated into
\begin{equation}\label{foncimplicite}
\partial_1 f_{(0,Y(0))} v + \partial_2 f_{(0,Y(0))} Y'(0) = (D\tau)_0 v.
\end{equation}
Lemma 6.6 in \cite{BK} ensures $\partial_2 f_{(0,Y(0))}$ is $\Lambda^2$-close to the identity. Besides, $\tau$ is an isometry for $\hat{g}$, 
so $(D\tau)_0$ is uniformly bounded. Finally, $\partial_1 f_{(0,Y(0))} v$ is the value at time $1$ of the Jacobi field $K$ along $s \mapsto s Y(0)$ corresponding to the geodesic variation
$$
H(t,s) \mapsto \Exp_{t v} s p(tv) Y(0).
$$ 
As the initial data for $K(s) = \partial_t H(0,s)$ are $K(0)=v$ and $K'(0)=0$, we obtain (corollary 6.3.8 of \cite{BK}) a bound on $K$ and thus on $\partial_1 f_{(0,Y(0))} v$. This yields a bound on $Y'(0)$ (thanks to \ref{foncimplicite}) and we are done.
\endproof

This lemma and the curvature decay lead to the estimate
$$
\abs{p_\gamma(t_2)^{-1} H_{\gamma(t_2)} p_\gamma(t_2) - p_\gamma(t_1)^{-1} H_{\gamma(t_1)} p_\gamma(t_1)} \leq C \int_{t_1}^\infty t^{-3} dt \leq C \, t_1^{-2}.
$$
Now recall the holonomy of the loops under consideration goes to the identity at infinity. Setting $t_1=:t$ and letting $t_2$ go to infinity, we find 
\begin{equation}\label{cvghol1}
\abs{p_\gamma(t)^{-1} H_{\gamma(t)} p_\gamma(t) - \id} \leq C \, t^{-2}.
\end{equation}

Since $M$ has zero (hence nonnegative) Ricci curvature and cubic volume growth, it follows from Cheeger Gromoll theorem that $M$ has only one end. Relying on faster than quadratic curvature decay, \cite{Kas} then ensures large spheres $S(o,t)$ are connected with intrinsic diameter bounded by $C s$. Thus every point $x$ in $S(o,t)$ is connected to $\gamma(t)$ by some curve $\beta$ with length at most $C t$ et and remaining outside $B(o,t/2)$. Holonomy comparison lemma \ref{compholo} yields:
\begin{equation}\label{cvghol2}
\abs{p_\beta^{-1} H_x p_\beta - H_{\gamma(t)} } \leq C t^{-2},
\end{equation} 
where $p_\beta$ is the parallel transportation along $\beta$ and $H_x$ is the holonomy endomorphism corresponding to a consistent orientation of $\sigma_x$. It follows that
\begin{equation}\label{cvghol3}
\abs{H_x - id} \leq C r(x)^{-2}.
\end{equation} 

So we have managed to improve our estimate on the holonomy of fundamental loops.

\begin{lem}
Let $(M^4,g)$ be a complete hyperk\"ahler manifold with
$$
\int_M \abs{\Rm}^2 r dvol < \infty
$$
and  
$$
\forall \, x \in M, \, \forall \, t \geq 1, \, A t^3 \leq \vol B(x,t) \leq B t^3
$$
($0 < A \leq B$). Then the holonomy $H_x$ of the fundamental loops $\sigma_x$ satisfies 
$$
\abs{H_x - \id} \leq C r(x)^{-2}. 
$$
\end{lem}


\subsection{Local Gromov-Hausdorff approximations.}

A first way to describe the local geometry consists in saying that, ``seen from far away'', it is close to a simpler geometry. We will show that the local geometry of the codimension $1$ collapsings that we are looking at is close to the Euclidean geometry in the immediately inferior dimension, for the Gromov-Hausdorff topology. 

\begin{rem}
Cubic curvature decay is important in the following lemma. A (not so) heuristic reason is the following. Forget geodesic loops and look at the metric in the exponential chart at $x$. If the curvature is bounded by $\Lambda^2$ in the area under consideration, comparison asserts
$$
\left( \frac{\sin \Lambda r}{\Lambda r} \right)^2 g_x \leq \exp_x^* g \leq \left( \frac{\sinh \Lambda r}{\Lambda r} \right)^2 g_x
$$
on a scale $r << \Lambda^{-1}$. The corresponding distances thus obey
$$
\frac{\sin \Lambda r}{\Lambda r} d_{g_x} \leq d_{\exp_x^* g} \leq \frac{\sinh \Lambda r}{\Lambda r} d_{g_x},
$$
hence
$$
\abs{d_{\exp_x^* g} - d_{g_x} } \leq C \Lambda^2 r^2 d_{g_x}.
$$
If we want to control the difference between these distances by some constant on the scale $r$, we therefore need a bound on $\Lambda^2 r^2 r$, which means the curvature ($\Lambda^2$) should be bounded by $r^{-3}$.  
\end{rem}

The following lemma ensures the elements of the fundamental group are almost translations.

\begin{lem}\label{taut}
Let $(M^4,g)$ be a complete hyperk\"ahler manifold with
$$
\int_M \abs{\Rm}^2 r dvol < \infty
$$
and  
$$
\forall \, x \in M, \, \forall \, t \geq 1, \, A t^3 \leq \vol B(x,t) \leq B t^3
$$
($0 < A \leq B$). Then there exists a compact set $K$ in $M$ and geometric constants $J$, $L$, $\kappa>0$ such that for every point $x$ in $M \backslash K$ and
every $\tau$ in $\Gamma(x,\kappa r(x))$, one has
$$
\Forall{w}{\hat{B}(0,\kappa r(x))} \abs{\tau (w)- t_{k v_x} (w)} \leq J
$$
where $v_x$ is a lift of the tip of $\sigma_x$ and $k$ is a natural number bounded by $L r(x)$. 
\end{lem}

\proof
Proposition \ref{pgf} asserts we can write $\tau = \tau_{v_x}^k$, where $v_x$ is a lift of a tip of $\sigma_x$ and $k$ is a natural number bounded by  $L r(x)$. Lemma \ref{holetpgf} ensures that for every $w$ in $\hat{B}(0,r(x)/4)$:
$$
\abs{\tau_{v_x}(w) - v_x - p_{v_x}^{-1}(w)} \leq C r(x)^{-3} \abs{v_x} \abs{w} (\abs{v_x} + \abs{w}).
$$
Thanks to cubic curvature decay, (\ref{cvghol3}) yields: 
$$
\abs{p_{v_x}^{-1}(w) - w} \leq C r(x)^{-2} \abs{w}.
$$
Combining these estimates, we obtain:
$$
\abs{\tau_{v_x}(w) - t_{v_x} (w)} = \abs{\tau_{v_x}(w) - v_x - w} \leq C r(x)^{-2} \abs{w}. 
$$
For every natural number $i \leq k$, we set $e_i = \tau_{v_x}^i - t_{v_x}^i$ and observe the formula
$$
e_{i+1} - e_i = e_1 \circ \tau_{v_x}^i.
$$
With
$$
\abs{\tau_{v_x}^i(w)} = d(\tau_{v_x}^i(w),0)
= d(\tau_{v_x}^{-i}(0),w) \leq \abs{\tau_{v_x}^{-i}(0)} + \abs{w},
$$
we find that for every $w$ in $\hat{B}(0,\kappa r(x))$:
$$
\abs{e_{i+1} (w) - e_i(w)} \leq C r(x)^{-2} \abs{\tau_{v_x}^i(w)} 
\leq C r(x)^{-1}.
$$  
By induction, it follows that $\abs{e_k (w)} \leq C k r(x)^{-1}$ and since $k \leq L r(x)$, we are led to: 
$$
\abs{\tau(w) - t_{k v_x} (w)} = \abs{\tau_{v_x}^k(w) - k v_x - w} = \abs{e_k (w)}
\leq C.
$$
\endproof

\begin{prop}[Gromov-Hausdorff approximation]\label{approxGH}
Let $(M^4,g)$ be a complete hyperk\"ahler manifold with
$$
\int_M \abs{\Rm}^2 r dvol < \infty
$$
and  
$$
\forall \, x \in M, \, \forall \, t \geq 1, \, A t^3 \leq \vol B(x,t) \leq B t^3
$$
($0 < A \leq B$). Then there exists a compact set $K$ in $M$ and geometric constants $I$, $\kappa>0$ such that every point $x$ in $M \backslash K$ has a neighborhood $\Omega$ whose Gromov-Hausdorff distance to the ball of radius $\kappa r(x)$ in $\Rl^{3}$ is bounded by $I$. 
\end{prop}

\proof
Choose a lift of $\sigma_x$ in $T_xM$ and denote by $v_x$ its tip. We call $H$ the hyperplane orthogonal to $v_x$ and write $v \mapsto v_H$ for the 
Euclidean orthogonal projection onto $H$ (for $g_x$). 

If $y$ is a point in $B(x,\kappa r(x)/2)$, we can define $h(y)$ as affine center of mass of the points $v_H$ obtained from lifts $v$ of $y$ in $\hat{B}(0,\kappa r(x)/2)$. This defines a map $h$ from $B(x,\kappa r(x)/2)$ to $H \cong \Rl^{3}$ (we endow $H$ of the Euclidean structure induced by $g_x=\abs{.}^2$). 

We consider the ball $B$ centered in $0$ and with radius $0.1 \kappa r(x)$ in $H$: $0.1 \kappa$ will be the $\kappa$ of the statement. Let us set $\Omega:= h^{-1} (B)$. We want to see that $\appli{h}{\Omega}{B}$ is the promised Gromov-Hausdorff approximation. We need to check that this map $h$ has $I$-dense image and that for all points $y$ and $z$ in $\Omega$:
$$
\abs{d(y,z) - \abs{h(y)-h(z)}} \leq I. 
$$
Firstly, since $v$ is in $B$, lemma \ref{taut} ensures that for every $\tau = \tau_{v_x}^k$ in $\Gamma(x,\kappa r(x))$, we have 
$$
\abs{\tau(v) - v - k v_x} \leq J 
$$
and thus, using Pythagore theorem,
$$
\abs{\tau(v)_H - v} \leq J. 
$$
Passing to the center of mass, we get
$$
\abs{h(\exp_x v) - v} \leq J. 
$$
If $d(v,H \backslash B) > J$, this proves $h(\exp_x v)$ belongs to $B$ and therefore $\exp_x v$ belongs to $\Omega$; as a result, $d(v,h(\Omega)) \leq J$. As $\set{v \in B \, / \, d(v,H \backslash B) > J}$ is $J$-dense in $B$, we have shown that $h(\Omega)$ is $2J$-dense in $B$. 

Secondly, consider two points $y$ and $z$ in $\Omega$. Lift them into $v$ and $w$ ($ \in B(x,\kappa r(x)/2)$) with $d(v,w) = d(y,z)$. As above, we get $\abs{h(y) - v_H} \leq J$ and $\abs{h(z) - w_H} \leq J$, hence
$$
\abs{\abs{h(y) - h(z)} - \abs{v_H - w_H}} \leq 2J.
$$
In particular, we obtain
$$
\abs{h(y) - h(z)} \leq \abs{v_H - w_H} + 2J \leq \abs{v - w} + 2J.
$$
Comparison yields
$$
\abs{v - w} \leq \frac{C r(x)^{-\frac{3}{2}} r(x)}{\sin C r(x)^{-\frac{3}{2}} r(x)} d(v,w) \leq \left(1 + C r(x)^{-1} \right) d(v,w)
$$
hence
$$
\abs{v - w} \leq d(v,w) + C r(x)^{-1} d(v,w) \leq d(v,w) + C.
$$
We deduce
$$
\abs{h(y) - h(z)} \leq d(v,w) + 2J + C = d(y,z) + 2J + C.
$$
Now, consider lifts $v'$ and $w'$ at minimal distance from $H$ and observe lemma \ref{trigo} yields:
$\abs{v'-v'_H} \leq C$ et $\abs{w'-w'_H} \leq C$. We deduce:
$$
\abs{v'-w'} \leq \abs{v'_H - w'_H} + C.
$$ 
The distance between $y$ and $z$ is nothing but the infimum of the distances between their lifts, so $d(y,z) \leq d(v',w')$. As above, comparison ensures:
$$
d(v',w') \leq \abs{v' - w'} + C.
$$
These three inequalities give altogether:
$$
d(y,z) \leq \abs{v'_H - w'_H} + C
$$
And since 
$$
\abs{\abs{h(y) - h(z)} - \abs{v'_H - w'_H}} \leq 2J,
$$
we arrive at 
$$
d(y,z) \leq \abs{h(y) - h(z)} + 2J + C.
$$
We have proved 
$$
\abs{d(y,z) - \abs{h(y) - h(z)}} \leq I,
$$
hence the result.
\endproof

The following step consists in regularizing local Gromov-Hausdorff approximations to obtain local fibrations which accurately describe the local geometry at infinity.  


\subsection{Local fibrations.}

The local Gromov-Hausdorff approximation that we built above has no reason to be regular. We will now smooth it into a fibration. The technical device  is simply a convolution, as in \cite{Fuk} and \cite{CFG}. We basically need theorem 2.6 in \cite{CFG}. The trouble is this general result will have to be refined, by using fully the cubic decay of the curvature and the symmetry properties of the special Gromov-Hausdorff approximation we smooth. This technique requires a control on the covariant derivatives of the curvature, but it is heartening to know that this is given for free on gravitational instantons (see theorem \ref{curvdec} in appendix A).

We say $f$ is a $C$-almost-Riemannian submersion if $f$ is a submersion such that for every horizontal vector $v$ (i.e. orthogonal to fibers),  
$$
e^{- C} \abs{v} \leq \abs{df_x(v)} \leq e^{C} \abs{v}.
$$  

\begin{prop}[Local fibrations]\label{fibrloc}
Let $(M^4,g)$ be a complete hyperk\"ahler manifold with
$$
\int_M \abs{\Rm}^2 r dvol < \infty
$$
and  
$$
\forall \, x \in M, \, \forall \, t \geq 1, \, A t^3 \leq \vol B(x,t) \leq B t^3
$$
($0 < A \leq B$). Then there exists a compact set $K$ in $M$ and geometric constants $\kappa>0$, $C >0$ such that for every point $x$ in $M \backslash K$, there is a circle fibration $\appli{f_x}{\Omega_x}{B_x}$ defined on a neighborhood $\Omega_x$ of $x$ and with values in the Euclidean ball $B_x$ with radius $\kappa r(x)$ in $\Rl^{3}$. Moreover,
\begin{itemize}
\item $f_x$ is a $C r(x)^{-1}$-almost-Riemannian submersion,
\item its fibers are submanifolds diffeomorphic to $\Sph^1$, with length pinched between $C^{-1}$ and $C$, 
\item $\abs{\nabla^2 f_x} \leq C r(x)^{-2}$,
\item $\forall i \geq 3, \abs{\nabla^i f_x} = \OO (r(x)^{1-i})$.
\end{itemize}
\end{prop}

\proof
In the proof of \ref{approxGH}, we introduced a function $h$ from the ball $B(x,\kappa r(x))$ to the hyperplane $H$, orthogonal to the tip $v_x$ of a lift of $\sigma_x$ in $T_x M$; this hyperplane $H$ is identified to the Euclidean space $\Rl^{3}$ through the metric induced by $g_x$. 

Let us choose a smooth nonincreasing function $\chi$ from $\Rl_+$ to $\Rl_+$, equal to $1$ on $[0,1/3]$ and $0$ beyond $2/3$. We also fix a scale $\epsilon := 0.1 \kappa r(x)$ and set $\chi_\epsilon(t) = \chi(2 t/\epsilon^2)$. Note the estimates:
\begin{equation}\label{estimchi}
\abs{\chi_\epsilon^{(k)}} \leq C_k \epsilon^{-2k}. 
\end{equation}

We consider the function defined on $B(x,\kappa r(x))$ by:
$$
f(y) := \frac{\int_{T_y M} h(\exp_y v) \chi_{\epsilon}(d(0,v)^2/2) dvol(v)}{\int_{T_y M} \chi_{\epsilon}(d(0,v)^2/2) dvol(v)}.
$$  
Here, $dvol$ and $d$ are taken with respect to $\exp_y^*g$. If $w$ is a lift of $y$ in $T_xM$, we can change variables thanks to the isometry
$$
\tau_w := \Exp_w \circ \left( T_w \exp_x \right)^{-1}
$$
between $(T_y M, \exp_y^*g)$ and $(T_x M, \exp_x^*g)$. For every point $v$ in $T_xM$, we introduce the function $\rho_v := \frac{d(v,.)^2}{2}$ and set $\hat{f} := f \circ \exp_x$, $\hat{h} := h \circ \exp_x$. We then get the formula: 
$$
f(y) = \hat{f}(w) = \frac{\int_{T_x M} \hat{h}(v) \chi_{\epsilon}(\rho_v(w)) dvol(v)}{\int_{T_x M} \chi_{\epsilon}(\rho_v(w)) dvol(v)}.
$$
The point is we can now work on a fixed Euclidean space, $(T_xM,g_x)$. The Riemannian measure $dvol$ can be compared to Lebesgue measure $dv$: on a scale $\epsilon$, if the curvature is bounded by $\Lambda^2$, we have 
$$
\left( \frac{\sin \Lambda \epsilon}{\Lambda \epsilon} \right)^4 dv \leq dvol \leq \left( \frac{\sinh \Lambda \epsilon}{\Lambda \epsilon} \right)^4 dv.
$$
Cubic curvature decay implies $\Lambda$ is of order $\epsilon^{-\frac{3}{2}}$, so that we find
\begin{equation}\label{compdvol}
- C \epsilon^{-1} dv \leq dvol - dv  \leq C \epsilon^{-1} dv.
\end{equation}
Distance comparison yields in the same way: 
$$
\abs{d(v,w) - \abs{v-w} } \leq C \Lambda^2 \epsilon^2 d(v,w) \leq C,
$$
hence
\begin{equation}\label{compdist2}
\abs{\rho_v(w)-\abs{v-w}^2/2} \leq C \epsilon.
\end{equation}
Eventually, the proof of \ref{approxGH} shows $\hat{h}$ is close to a Euclidean projection onto $H$:
\begin{equation}\label{compproj}
\abs{\hat{h}(v) - v_H} \leq C.
\end{equation}
We can write
\begin{eqnarray*}
\int \hat{h}(v) \chi_{\epsilon}(\rho_v(w)) dvol(v)
&=& \int \hat{h}(v) \chi_{\epsilon}(\rho_v(w)) (dvol(v) -dv) \\
&+& \int \hat{h}(v) \left(\chi_{\epsilon}(\rho_v(w)) -  \chi_{\epsilon}(\abs{v-w}^2/2)\right) dv \\
&+& \int (\hat{h}(v) -v_H) \chi_{\epsilon}(\abs{v-w}^2/2) dv \\
&+& \int v_H \chi_{\epsilon}(\abs{v-w}^2/2) dv.
\end{eqnarray*}
The support of $v \mapsto \chi_{\epsilon}(\rho_v(w))$ is included in a ball whose radius is of order $\epsilon$: $\hat{h}$ will therefore take its values in a ball with radius of order $\epsilon$. With (\ref{compdvol}), we can then bound the first term of the right-hand side by $C \epsilon \cdot \epsilon^{-1} \cdot \epsilon^4 = C \epsilon^4$. With (\ref{estimchi}) and (\ref{compdist2}), we bound the second term by $C \epsilon \cdot \epsilon^{-2} \cdot \epsilon \cdot \epsilon^4 = C \epsilon^4$. Eventually, (\ref{compproj}) controls the third term by $C \epsilon^4$. We get:
$$
\int \hat{h}(v) \chi_{\epsilon}(\rho_v(w)) dvol(v) = \int v_H    \chi_{\epsilon}(\abs{v-w}^2/2) dv + \OO(\epsilon^4),
$$
where $\OO(\epsilon^4)$ stands for an error term of magnitude $\epsilon^4$. 

Thanks to (\ref{compdvol}), (\ref{estimchi}) and (\ref{compdist2}), we obtain in the same way:
$$
\int_{T_x M} \chi_{\epsilon}(\rho_v(w)) dvol(v) = \int \chi_{\epsilon}(\abs{v-w}^2/2) dv + \OO(\epsilon^{3}).
$$
Observing
$$
\int v_H  \chi_{\epsilon}(\abs{v-w}^2/2) dv = \OO(\epsilon^{5})
$$
and
$$
C^{-1} \epsilon^4 \leq \int \chi_{\epsilon}(\abs{v-w}^2/2) dv \leq C \epsilon^4,
$$
we deduce 
$$
\hat{f}(w) = \frac{\int v_H \chi_{\epsilon}(\abs{v-w}^2/2) dv}{\int \chi_{\epsilon}(\abs{v-w}^2/2) dv} +  \OO(1).
$$
The change of variables $z=v-w$ yields:
$$
\hat{f}(w) - w_H = \underbrace{\frac{\int z_H \chi_{\epsilon}(\abs{z}^2/2) dz}{\int \chi_{\epsilon}(\abs{z}^2/2) dz}}_{=0 \text{ by parity}} + \OO(1),
$$
hence
\begin{equation}\label{compfib}
\hat{f}(w) = w_H + \OO(1).
\end{equation}

The differential of $\hat{f}$ reads
$$
d\hat{f}_w = \frac{\int ( \hat{h}(v) - \hat{f}(w) ) \chi'_{\epsilon}(\rho_v(w)) (d\rho_v)_w dvol(v)}{\int \chi_{\epsilon}(\rho_v(w)) dvol(v)}.
$$
The same kind of approximations, based on (\ref{compvol}), (\ref{estimchi}), (\ref{compdist2}), (\ref{compderdist}), (\ref{compproj}) and (\ref{compfib}) imply 
\begin{equation}\label{compderfib}
d\hat{f}_w = - \frac{\int z_H \chi'_{\epsilon}(\abs{z}^2/2) (z,.) dz}{\int \chi_{\epsilon}(\abs{z}^2/2) dz} + \OO(\epsilon^{-1}).
\end{equation}

Let us choose an orthonormal basis $(e_1,\dots,e_4)$ of $T_xM$, with $e_4 \bot H$. If $i \not= j$, parity shows 
$$
\int z_i \chi'_{\epsilon}(\abs{z}^2/2) z_j dz = 0.
$$
On the contrary, an integration by parts ensures that for every $\alpha \geq 0$:
$$
\int_{-\infty}^\infty z_i^2 \chi'_{\epsilon}(z_i^2/2 + \alpha) dz_i = - \int_{-\infty}^\infty \chi_{\epsilon}(z_i^2/2 + \alpha) dz_i,
$$
so that 
$$
-\int z_i^2 \chi'_{\epsilon}(\abs{z}^2/2) dz = \int \chi_{\epsilon}(\abs{z}^2/2) dz.
$$
This means precisely:
$$
- \frac{\int z_H \chi'_{\epsilon}(\abs{z}^2/2) (z,.) dz}{\int \chi_{\epsilon}(\abs{z}^2/2) dz} = \sum_{i=1}^{3} e_i \otimes (e_i,.).
$$
And one can recognize the Euclidean projection onto $H$. We deduce 
$$
d\hat{f}_w = \sum_{i=1}^{3} e_i \otimes (e_i,.)  + \OO(\epsilon^{-1}) 
$$
With (\ref{compmet}), this proves $\hat{f}$ is a $C \epsilon^{-1}$-almost-Riemannian submersion. Since $\exp$ is a local isometry, $f$ 
is also a $C \epsilon^{-1}$-almost-Riemannian submersion.  

The Hessian reads:
\begin{eqnarray*}
\nabla^2 \hat{f}_w &=& \frac{\int ( \hat{h}(v) - \hat{f}(w) ) \left( \chi''_{\epsilon}(\rho_v(w)) (d\rho_v)_w\otimes (d\rho_v)_w + \chi'_{\epsilon}(\rho_v(w)) (\nabla^2 \rho_v)_w \right) dvol(v)}{\int \chi_{\epsilon}(\rho_v(w)) dvol(v)} \\
&-& 2 d\hat{f}_w \otimes \frac{\int \chi'_{\epsilon}(\rho_v(w)) (d\rho_v)_w) dvol(v)}{\int \chi_{\epsilon}(\rho_v(w)) dvol(v)}.
\end{eqnarray*}

Again, with (\ref{compvol}), (\ref{estimchi}), (\ref{compdist2}), (\ref{compderdist}), (\ref{compproj}) and (\ref{compfib}), we arrive at 
\begin{eqnarray*}
\nabla^2 \hat{f}_w &=& \frac{\int  z_H  \left( \chi''_{\epsilon}(\abs{z}^2/2) (z,.)\otimes (z,.) + \chi'_{\epsilon}(\abs{z}^2/2) (.,.) \right) dz}{\int \chi_{\epsilon}(\abs{z}^2/2) dz} \\
&-& 2 d\hat{f}_w \otimes \frac{\int \chi'_{\epsilon}(\abs{z}^2/2) (z,.) dz}{\int \chi_{\epsilon}(\abs{z}^2/2) dz} + \OO(\epsilon^{-2}).
\end{eqnarray*}
To begin with, parity ensures
$$
\int \chi'_{\epsilon}(\abs{z}^2/2) (z,.) dz = 0
\quad \text{and}
\quad
\int  z_H  \chi'_{\epsilon}(\abs{z}^2/2) (.,.) dz = 0.
$$
The $i$th component of the integral
$$
\int  z_H  \chi''_{\epsilon}(\abs{z}^2/2)  (z,.)\otimes (z,.)
$$
can be written as a sum of terms 
$$
\left( \int  z_i z_j z_k  \chi''_{\epsilon}(\abs{z}^2/2) dz_1 \dots dz_4 \right) (e_j,.) \otimes (e_k,.)
$$
which vanish for a parity reason. Therefore:
$$
\nabla^2 \hat{f}_w = \OO(\epsilon^{-2}).
$$

The proof of theorem 2.6 in \cite{CFG} yields the remaining properties of $f_x:=f$. Essentially, $f$ is a fibration because it is $C^1$-close to a fibration. The connexity of the fibers follows from the bound on the Hessian of $f$. The length of the fibers is controlled by the assumption on the volume growth (since $f$ is a almost-Riemannian submersion). 
\endproof

We will need to relate neighboring fibrations (this somewhat corresponds to proposition $5.6$ in \cite{CFG}). 

\begin{lem}[Closeness of local fibrations I]\label{fibrproches}
The setting is the same as in proposition \ref{fibrloc}. Given two points $x$ and $x'$ in $M \backslash K$, with $ d(x,x') \leq \kappa r(x)$, if $\Omega_{x,x'}= \Omega_x \cap \Omega_{x'}$ (notations in \ref{fibrloc}), then there is a $C r(x)^{-1}$-almost-isometry $\phi_{x,x'}$ between 
$f_{x'} (\Omega_{x,x'})$ and $f_{x} (\Omega_{x,x'})$, for which moreover  
\begin{itemize}
\item $\abs{f_{x} - \phi_{x,x'} \circ f_{x'}} \leq C$,
\item $\abs{Df_{x} - D\phi_{x,x'} \circ Df_{x'}} \leq C r(x)^{-1}$,
\item $\abs{D^2\phi_{x,x'}} \leq C r(x)^{-2}$,
\item $\forall \, i \geq 3, \abs{D^i \phi_{x,x'} } = \OO (r(x)^{1-i})$.
\end{itemize}
\end{lem}

\proof
We use the same notations as in the previous proof, adding subscripts to precise the point under consideration, and we work in $T_xM$. Choose a lift $u$ of $y$ at minimal distance from $o$ and set $\tau_u := \Exp_u \circ (T_u \exp_x)^{-1}$ the corresponding isometry (between large balls in  $T_{x'}M$ and $T_xM$). We consider the map
$$
\phi_{x,x'} := f_x \circ \exp_{x'} \restric{f_{x'} (\Omega_{x,x'})}.
$$
In order to bring everything back into $T_x M$, we write  
$$
\phi_{x,x'} \circ f_{x'} \circ \exp_x = f_x \circ \exp_{x'} \circ f_{x'} \circ \exp_x.
$$
The relation $\exp_x \circ \tau_u = \exp_{x'}$ leads to the reformulation  
$$
\phi_{x,x'} \circ f_{x'} \circ \exp_x = f_x \circ \exp_x \circ \tau_u \circ f_{x'} \circ \exp_{x'} \circ \tau_u^{-1}
$$
that is
\begin{equation}\label{relphi}
\phi_{x,x'} \circ f_{x'} \circ \exp_x = \hat{f}_x \circ \tilde{f}_{x'}
\end{equation}
with $\hat{f}_x = f_x \circ \exp_x$ and $\tilde{f}_{x'} = \tau_u \circ f_{x'} \circ \exp_{x'} \circ \tau_u^{-1}$. We need to understand this lattest map.

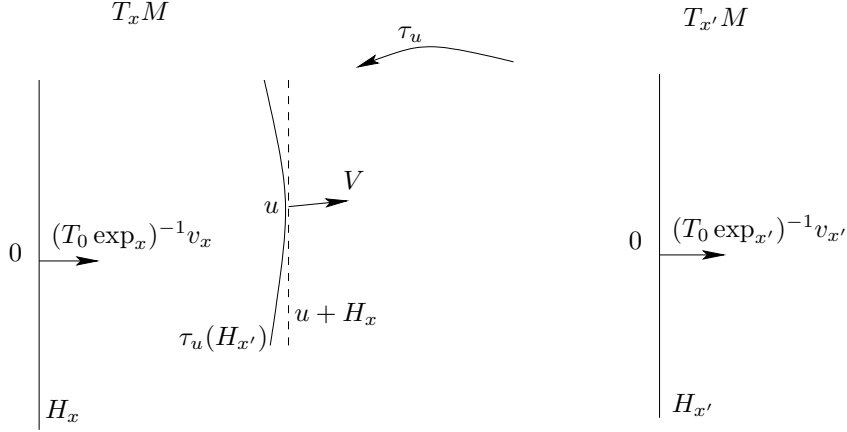
\begin{figure}[htb!]
\begin{center}
\input{compfibr.pstex_t}
\end{center}
\label{compfibr}
\end{figure}

Since $\tau_u$ is an isometry between the metrics $\exp_{x'}^* g$ and $\exp_{x}^* g$ and since $H_{x'}$ is the union of all the geodesics starting from $0$ and with a unit speed orthogonal to $(T_0\exp_{x'})^{-1}(v_{x'})$, $\tau_u(H_{x'})$ is the hypersurface generated by the geodesics starting from  $u$ with a unit speed orthogonal to $V := (d\tau_u)_0 \circ (T_0\exp_{x'})^{-1}(v_{x'})$. $v_{x'}$ is by definition one of the lifts of $x'$ by $\exp_{x'}$ which are not $0$ but at minimal distance from $0$ (in $T_{x'} M$). So $\tau_u (v_x')$ is one of the two lifts of $x'$ by $\exp_{x}$ which are not $\tau_u(0)=u$ but at minimal distance from $\tau_u(0)=u$ (in $T_{x} M$). We have seen in lemma \ref{reglacet} that such a point $\tau_u (v_x')$ is $\tau_{v_x}(u)$ or $\tau_{v_x}^{-1}(u)$. To fix ideas, assume we are in the first case: $\tau_u (v_x') = \tau_{v_x}(u)$. 

The exponential map of $T_{x'} M$ (at $0$) maps $(T_0\exp_{x'})^{-1}(v_{x'})$ to $v_{x'}$, so $V = (d\tau_u)_0 \circ (T_0\exp_{x'})^{-1}(v_{x'})$ is 
the vector which is mapped by the exponential map of $T_xM$ (at $\tau_u(0)=u$) to $\tau_u (v_x') = \tau_{v_x}(u)$: $\Exp_ u V = \tau_{v_x}(u)$. Consider the geodesic $\gamma(t) := \Exp_u t V$. Taylor formula
$$
\gamma(1) - \gamma(0) - \dot{\gamma}(0) = \int_0^1 (1-t) \ddot{\gamma}(t) dt
$$ 
and the estimate $\abs{\ddot{\gamma}} \leq C r(x)^{-2} \abs{V}^2 \leq C r(x)^{-2}$, steming from lemma \ref{kaul} and the bound on the injectivity radius (\ref{majinj}), together imply
$$
\abs{\tau_{v_x} (u) - u - V} \leq C r(x)^{-2}.
$$
With the estimate
$$
\abs{\tau_{v_x}(u) - u - v_x} \leq C r(x)^{-1},
$$
we deduce 
\begin{equation}\label{estimvect}
\abs{V - v_x} \leq C r(x)^{-1}.
\end{equation}
The angle between vectors $V$ and $v_x$ is thus bounded by $C r(x)^{-1}$, so that, with $\hat{\U} := \hat{B}(0,\kappa r(x)) \cap \hat{B}(u,\kappa r(x'))$, the affine hyperplanes pieces $(u + V^\bot) \cap \hat{\U}$ and $(u + v_x^\bot) \cap \hat{\U}$ remain at bounded distance.

Considering the geodesic $\gamma(t) = \Exp_u tW$, with $W \bot V$ and $\abs{W} \leq C r(x)$, we obtain in the same way (thanks to lemma \ref{kaul}): 
$$
\abs{\Exp_u W - u - W} \leq C r(x)^{-2} r(x)^2 = C.
$$
This means the affine hyperplane piece $(u + V^\bot) \cap \hat{\U}$ and the hypersurface piece $\tau_u(B_{x'}) \cap \hat{\U} = \Exp_u V^\bot \cap \hat{\U}$ remain at bounded distance. 

And we conclude $\tau_u(B_{x'})\cap \hat{\U}$ and $(u + v_x^\bot) \cap \hat{\U}$ remain $C$-close, namely the map $\Psi$ defined from $\tau_u(B_{x'}) \cap \hat{\U}$ to $(u + v_x^\bot) \cap \hat{\U}$ by
$$
\psi : \, \Exp_u W \, \mapsto \, u + W
$$ 
is $C$-close to the identity.

In the preceding proof, we saw that $f_{x'} \circ \exp_{x'}$ was $C$-close to the orthogonal projection (for $g_{x'}$) onto $H_{x'}$. Now, $\tau_u$ is an isometry between the metrics $\exp_{x'}^* g$ and $\exp_{x}^* g$, which are respectively $C r(x)^{-1}$-close to $g_{x'}$ and $g_x$. Thus for every point $w$ in the area under consideration, $\abs{\tilde{f}_{x'}(w)-w}$ is $C$-close to the distance (for $g_{x}$) between $w$ and $\tau_u(H_{x'})$, so that $\abs{\psi \circ \tilde{f}_{x'}(w)-w}$ is $C$-close to the distance (for $g_{x}$) between $w$ and $(u + v_x^\bot)$: $\psi \circ \tilde{f}_{x'}$ and thus $\tilde{f}_{x'}$ are $C$-close to the orthogonal projection onto $(u + v_x^\bot)$, which is nothing but the composition of the orthogonal projection onto $H_x$ and of the translation with vector $u-u_{H_x}$:
$$
\abs{\tilde{f}_{x'} (w) - w_{H_x} - (u-u_{H_x}) } \leq C.
$$
We deduce 
$$
\abs{\tilde{f}_{x'} (w) - \hat{f}_x(w) - (u-u_{H_x}) } \leq C
$$
and, composing with $\hat{f}_x$, we find
$$
\abs{\hat{f}_x \circ \tilde{f}_{x'} (w) - \hat{f}_x (w)} \leq C.
$$
Recalling formula (\ref{relphi}), we obtain
$$
\abs{\phi_{x,x'} \circ f_{x'} \circ \exp_x - f_x \circ \exp_x} \leq C,
$$
and, with the surjectivity of $\exp_x$, this yields
$$
\abs{\phi_{x,x'} \circ f_{x'} - f_x} \leq C.
$$
Relation (\ref{relphi}) also implies
\begin{equation}\label{relderphi}
D (\phi_{x,x'} \circ f_{x'} \exp_x) = D\hat{f}_x \circ D\tilde{f}_{x'}.
\end{equation}
Let $z$ be a point in $\hat{\U}$ and set $z'= \tau_u^{-1} (z) \in T_{x'}M$. The preceding proof has shown that $D_z\hat{f}_x$ is $C r(x)^{-1}$-close to the orthogonal projection in the direction of $H_x$. In the same way, $D_{z'} (f_{x'}\exp_{x'})$ is $C r(x)^{-1}$-close to the orthogonal projection in the direction of $H_{x'}$, i.e. in the direction orthogonal to $v_{x'}$. Conjugating by $D\tau_u$, we find that $D_z \tilde{f}_{x'}$ is $C r(x)^{-1}$ close to the projection in the direction orthogonal to $D_z \tau_u (v_{x'})$. 

Let $Z'$ be the initial speed of the geodesic connecting $z'$ to $\tau_{v_{x'}} (z')$ in unit time. The argument leading to (\ref{estimvect}) yields
$$
\abs{Z' - v_{x'}} \leq C r(x)^{-1}.  
$$
If we set $Z:= D_z \tau_u Z'$, we thus have
$$
\abs{Z - D_z \tau_u (v_{x'})} \leq C r(x)^{-1}.  
$$
Now $Z$ is the initial speed of the geodesic connecting $z$ to $\tau_{v_x} (z)$ (or $\tau_{v_x}^{-1} (z)$) in unit time. So again:
$$
\abs{Z - v_x} \leq C r(x)^{-1},
$$
so that 
$$
\abs{v_x - D_z \tau_u (v_{x'})} \leq C r(x)^{-1}.  
$$
Finally, $D_z \tilde{f}_{x'}$ is $C r(x)^{-1}$-close to the projection in the direction of the hyperplane $H_x$, orthogonal to $v_x$: 
$$
\abs{D (\phi_{x,x'} \circ f_{x'} \circ \exp_x)- D\hat{f}_x } \leq C r(x)^{-1},
$$
hence
\begin{equation}\label{estimdphi}
\abs{D \phi_{x,x'} \circ D f_{x'} - D f_x } \leq C r(x)^{-1}.
\end{equation}

Let $W$ be a vector tangent to $f_{x'} (\Omega_{x,x'})$ and let $W'$ be its horizontal lift for $f_{x'}$: $D f_{x'} W' = W$. As $D \tilde{f}_{x'}$ and $D\hat{f}_x$ are $C r(x)^{-1}$-close, an horizontal vector for $f_{x'}$ is $C r(x)^{-1}$-close to a horizontal vector for $f_x$. And since $f_x$ and $f_{x'}$ are $C r(x)^{-1}$-almost-Riemannian submersions, we get
$$
\abs{ \abs{D f_x (W')} - \abs{W'} } \leq C r(x)^{-1} \abs{W'}
$$
and 
$$
\abs{ \abs{W} - \abs{W'} } \leq C r(x)^{-1} \abs{W'}.
$$
Writing 
\begin{eqnarray*}
& &\abs{\abs{D \phi_{x,x'} W} - \abs{W}} \\
&\leq& \abs{\abs{D \phi_{x,x'} (D f_{x'} W')} - \abs{D f_x W'}} 
+ \abs{\abs{D f_x W'}-\abs{W'}} + \abs{\abs{W'}-\abs{W}} \\
&\leq& \abs{D \phi_{x,x'} (D f_{x'} W') - D f_x W'} 
+ \abs{\abs{D f_x W'}-\abs{W'}} + \abs{\abs{W'}-\abs{W}}
\end{eqnarray*}
and using (\ref{estimdphi}), we obtain
$$
\abs{\abs{D \phi_{x,x'} (W)} - \abs{W}} \leq C r(x)^{-1} \abs{W},
$$
which proves $\phi_{x,x'}$ is a $C r(x)^{-1}$-quasi-isometry. 

Higher order estimates stem from those on $f_x$ and $f_{x'}$, thanks to formula (\ref{relphi}): the bounds on the curvature covariant derivatives ensure $\exp_{x}^* g$ (resp. $\exp_{x'}^* g)$) is close to the flat $g_{x}$ (resp. $g_{x'}$) in $C^\infty$ topology, so that the estimates for one or the other are equivalent; thus $\exp_{x}$, $\exp_{x'}$ and $\tau_u$ can be treated like isometries.
\endproof

We will also need the following lemma. Indeed, it stems from the previous one. 

\begin{lem}[Local fibration closeness II]\label{fibrproches2}
The setting is the same as in lemma \ref{fibrproches}. We consider three points $x$, $x'$ and $x''$ in $M \backslash K$, whose respective distances are bounded by $\kappa r(x)$. Then, wherever it makes sense, we have 
\begin{itemize}
\item $\abs{\phi_{x,x''} - \phi_{x,x'} \circ \phi_{x',x''}} \leq C$,
\item $\abs{D \phi_{x,x''} - D\phi_{x,x'} \circ D\phi_{x',x''}} \leq C r(x)^{-1}.$
\end{itemize}
\end{lem}

\proof
On the intersection of $\Omega_x$, $\Omega_{x'}$ and $\Omega_{x''}$, we can write 
$$
\abs{f_{x} - \phi_{x,x'} \circ f_{x'}} \leq C \quad \text{ et } \quad \abs{f_{x'} - \phi_{x',x''} \circ f_{x''}} \leq C.
$$ 
Since $\phi_{x,x'}$ is a quasi-isometry, it follows that:
$$
\abs{f_{x} - \phi_{x,x'} \circ \phi_{x',x''} \circ f_{x''}} 
\leq \abs{ f_{x} - \phi_{x,x'} \circ f_{x'}} + \abs{\phi_{x,x'} \circ ( f_{x'} -  \phi_{x',x''} \circ f_{x''})}
\leq C. 
$$ 
Using the estimate 
$$
\abs{f_{x} - \phi_{x,x''} \circ f_{x''}} \leq C,
$$ 
we obtain by triangle inequality:
$$
\abs{(\phi_{x,x''} - \phi_{x,x'} \circ \phi_{x',x''}) \circ f_{x''}} \leq C. 
$$
From the surjectivity of $f_{x''}$, we see that:
$$
f_{x''}(\Omega_{x,x''} \cap f_{x''}(\Omega_{x',x''} \cap \phi_{x',x''}^{-1} f_{x'} (\Omega_{x,x'}).
$$
Since $f_{x''}$ is a submersion, the same argument applies to the differentials.
\endproof


\subsection{Local fibration gluing.}

Now, we need to adjust the local fibrations so as to make them compatible. The technical device is essentially the same as in \cite{CFG}. The following lemma will be widely used in this process. 

\begin{lem}[Local fibration adjustment I]\label{ajustement}
The setting is that of lemma \ref{fibrproches}. Given two points $x$ and $x'$ in $M \backslash K$ with $\alpha r(x) \leq d(x,x') \leq \beta r(x)$ for some real numbers $0<\alpha<\beta<1$. We assume that on $B(x,\gamma r(x))$ and $B(x', \gamma r(x'))$, some fibrations $f_x$ and $f_{x'}$ as in \ref{fibrloc} are defined, that $B(x,\delta r(x))$ and $B(x', \delta r(x'))$ have nonempty intersection, with $0<\delta<\gamma$, and that a map $\phi_{x',x}$ as in \ref{fibrproches} is defined. We can then build a fibration $\tilde{f}_{x'}$ on $B(x', \delta r(x'))$, with the same properties as $f_{x'}$, plus:
$$
\tilde{f}_{x'} = \phi_{x',x} \circ f_{x}
$$
on $B(x,\delta r(x)) \cap B(x', \delta r(x'))$. Moreover, this new fibration coincides with the old $f_{x'}$ on $B(x,\gamma r(x))$ and wherever we already had $f_{x'} = \phi_{x',x} \circ f_{x}$. \end{lem}

\proof
We set 
$$
\tilde{f}_{x'}(y) = \lambda(y) \phi_{x',x} (f_{x}(y)) + (1-\lambda(y)) f_{x'}(y)
$$
with 
$$
\lambda(y) = \theta \left( \frac{f_{x}(y)}{r(x)} \right)
$$
where $\appli{\theta}{\Rl^{3}}{[0,1]}$ is a truncature function equal to $1$ on the ball centered in $0$ and with radius $\delta$, equal to $0$ outside the ball centered in $0$ and with radius $\gamma$. Using the bounds on $f_x$, we find $\abs{\nabla^k \lambda} \leq C_k r(x)^{-k}$,
and the announced estimates can be obtained by differentiating the equation 
$$
\tilde{f}_{x'}(y) - f_{x'}(y) = \lambda(y) \left( \phi_{x',x} \circ f_{x}(y)) -  f_{x'}(y) \right).
$$
\endproof

\begin{lem}[Local fibration adjustment II]\label{ajustement2}
The setting is that of lemma \ref{fibrproches2}. Given three points $x$,$x'$ and $x''$ in $M \backslash K$ with $\alpha r(x) \leq d(x,x'), \, d(x',x''), \, d(x,x'') \leq \beta r(x)$ for some real numbers $0<\alpha<\beta<1$. We assume that on $B(x,\gamma r(x))$, $B(x', \gamma r(x'))$ and $B(x'', \gamma r(x''))$, some fibrations $f_x$, $f_{x'}$ and $f_{x''}$ as in \ref{fibrloc} are defined, that the intersection of $B(x,\delta r(x))$, $B(x', \delta r(x'))$ and $B(x'', \delta r(x''))$ is nonempty for some $0<\delta<\gamma$ and that maps $\phi_{x',x}$, $\phi_{x,x''}$ and $\phi_{x',x''}$ as in \ref{fibrproches2} are defined. We can then build a new diffeomorphism $\tilde{\phi}_{x',x''}$, with the same properties as $\phi_{x',x''}$, plus:
$$
\tilde{\phi}_{x',x''} = \phi_{x',x} \circ \phi_{x,x''}
$$
on $f_{x''} (B(x,\delta r(x)) \cap B(x', \delta r(x'))\cap B(x'', \delta r(x''))$. Moreover, this new diffeomorphism coincides with $\phi_{x',x''}$ on $B(x'',\gamma r(x''))$ and wherever we already had $\phi_{x',x''} = \phi_{x',x} \circ \phi_{x,x''}$. 
\end{lem}

\proof
We simply set 
$$
\tilde{\phi}_{x',x''}(v) = \lambda(v) \phi_{x',x} \circ \phi_{x,x''} (v) + (1-\lambda(v)) \phi_{x',x''}(v)
$$
with 
$$
\lambda(v) = \theta \left( \frac{\abs{v}^2}{r(x)^2} \right)
$$
where $\theta$ is the same function as in the previous proof. 
\endproof

\begin{thm}[Global fibration]\label{fibration}
Let $(M^4,g)$ be a complete hyperk\"ahler manifold with
$$
\int_M \abs{\Rm}^2 r dvol < \infty
$$
and  
$$
\forall \, x \in M, \, \forall \, t \geq 1, \, A t^3 \leq \vol B(x,t) \leq B t^3
$$
($0 < A \leq B$). 
Then there exists a compact set $K$ in $M$ such that $M \backslash K$ is endowed with a smooth circle fibration $\pi$ over a smooth open manifold $X$.
Besides, there is a geometric positive constant $C$ such that fibers have length pinched between $C^{-1}$ and $C$ and second fundamental form bounded by $C r^{-2}$. 
\end{thm}

\begin{rem}
The proof will show that for any point $x$ in $M \backslash K$, there is a diffeomorphism $\psi_x$ between a neighborhood of $\pi(x)$ in $X$ and a ball in $\Rl^{3}$ such that $\psi_x \circ \pi$ is a fibration satisfying estimates as in proposition \ref{fibrloc}.
\end{rem}

\proof
We take a maximal set of points $x_i$, $i \in I$, such that for all indices $i \not= j$,  $d(x_i,x_j) \geq \kappa r(x_i)/8$. This provides a uniformly locally finite covering of $M$ by the balls $B(x_i,\kappa r(x_i)/2)$. For every index $i$, we let $f_i$ be the local fibration given by \ref{fibrloc}. We will work with the minimal saturated (for $f_i$) sets $\Omega_i(\alpha)$ containing the balls $B(x_i, \alpha r(x_i))$, where $\alpha$ is a parameter inferior to $\kappa$. As in \cite{CFG}, we divide $I$ into packs $S_1$, ..., $S_N$ such that any two distinct points $x_i$, $x_j$ whose indices are in the same pack are far from each other: 
$$
\Exists{a}{[1,N]} \set{i,j} \subset S_a \Rightarrow d(x_i,x_j) \geq 100 \kappa \min(r(x_i),r(x_j)).
$$ 
In particular, $\Omega_i(\alpha)$ and $\Omega_j(\alpha)$ have empty intersection if $i$ and $j$ are in different packs; in this case, if the number of the pack of $i$ is greater than for $j$, one denotes by $\phi_{i,j}$ the diffeomorphism given by \ref{fibrproches} and by $\phi_{j,i}$ its inverse. 

In order to improve the approximations $f_i \approx \phi_{i,j} f_j$ into equalities $f_i = \phi_{i,j} f_j$, we set up an adjustment campaign in the following way. The idea consists in giving priority to packs with small number. To do so, given an area where several fibrations are defined, we will modify them so that they all fit with the fibration with smallest number among them. The order of implementation is important. We will distinguish several stages, indexed by subsets $\A:=\set{a_1 < \cdots < a_k}$ of $[1,N]$. We implement these $2^N$ stages by increasing order of $a_1$, then decreasing order of $k$, then increasing order of $a_2$, then increasing order of $a_3$, etc. To rephrase it, we have 
$$
\set{a_1 < \cdots < a_k} \prec \set{b_1 < \cdots < b_l} 
$$
if one of these exclusive conditions is realized:
\begin{itemize}
\item $a_1 < b_1$ ;
\item $a_1=b_1$  and  $k > l$ ;
\item $a_i=b_i$ for $i \leq i_0$  and  $k = l$  and  $a_{i_0} < b_{i_0}$. 
\end{itemize}
We denote by $m_\A$ the rank of $\A$ in this order and set
$$
\alpha_m := \kappa \cdot \left( \frac{1}{2} \right)^{\frac{m}{2^N}}.
$$
Along the campaign, the fibration domains $\Omega_i(\alpha)$ will be shrinked: $\alpha_{m_\A}$ will be the domain size at stage $\A$.  

At stage $\A:=\set{a_1 < \cdots < a_k}$, we consider all elements $\I=(i_1,\cdots,i_k)$ of $S_{a_1} \times \cdots \times S_{a_k}$: to each such element corresponds one step. At step $\I$, we are interested in $\Omega_\I:=\Omega_{i_1}(\alpha_{m_\A+1}) \cap \cdots \cap \Omega_{i_k}(\alpha_{m_\A+1})$. One should notice that our choice of packing ensures all the intersections $\Omega_{i_1}(\alpha_{m_\A}) \cap \cdots \cap \Omega_{i_k}(\alpha_{m_\A})$ treated at the same stage are away from each other, so that the following modifications are independent (during the stage). Essentially, the fibration $f_{i_1}$ will overrule its neighbour on $\Omega_\I$. Given $2 \leq p \leq k$, we build $\tilde{f}_{i_p}$ on $\Omega_{i_p}(\alpha_{m_\A+1})$, from $f_{i_1}$ and $f_{i_p}$, as in \ref{ajustement}, so as to obtain 
\begin{itemize}
\item $\tilde{f}_{i_p} = \phi_{i_p,i_1} f_{i_1} \text{ sur } \Omega_{i_p}(\alpha_{m_\A+1}) \cap \Omega_{i_1}(\alpha_{m_\A+1})$, 
\item $\tilde{f}_{i_p} = f_{i_p} \text{ sur } \Omega_{i_p}(\alpha_{m_\A+1}) \backslash \Omega_{i_1}(\alpha_{m_\A})$. 
\end{itemize}
We also build, for $2\leq p < q \leq k$, $\tilde{\phi}_{i_p,i_q}$ on $\tilde{f}_{i_q} ( \Omega_{i_p}(\alpha_{m_\A+1}) \cap \Omega_{i_q}(\alpha_{m_\A+1}) )$ from $\phi_{i_p,i_1} \phi_{i_1,i_q}$ and $\phi_{i_p,i_q}$, as in \ref{ajustement2}, so that 
\begin{itemize}
\item $\tilde{\phi}_{i_p,i_q} = \phi_{i_p,i_1} \phi_{i_1,i_q} \text{ sur } \tilde{f}_{i_q} ( \Omega_{i_p}(\alpha_{m_\A+1}) \cap \Omega_{i_q}(\alpha_{m_\A+1}) \cap \Omega_{i_1}(\alpha_{m_\A})  )$, 
\item $\tilde{\phi}_{i_p,i_q} = \phi_{i_p,i_q} \text{ sur } \tilde{f}_{i_q} ( \Omega_{i_p}(\alpha_{m_\A+1}) \cap \Omega_{i_q}(\alpha_{m_\A+1}) \backslash \Omega_{i_1}(\alpha_{m_\A})  )$. 
\end{itemize}
After this, we can add that wherever it makes sense, we have for every $\set{p,q} \subset [2,k]$:
$$
\tilde{\phi}_{i_q,i_p} \tilde{f}_{i_p} = \phi_{i_q,i_1} \phi_{i_1,i_p} \phi_{i_p,i_1} f_{i_1} = \phi_{i_q,i_1} f_{i_1} = \tilde{f}_{i_q}.
$$
Now forget the tildes. We have just ensured that on $\Omega_\I$, for all relevant indices $i, j$, one has $f_i = \phi_{i,j} f_j$.

We proceed, independently, for all possible $\I$ at this stage, then we go on with the next stage, following the chosen order.

At the moment we pass from a stage $\set{a_1 <\cdots }$ to a stage $\set{b_1 <\cdots }$, with $a_1 \not= b_1$, we can notice the fibrations $f_i$ and the diffeomorphisms $\phi_{i,j}$ are definitively fixed on the sets with number in the pack $S_{a_1}$: indeed, the device of \ref{ajustement} and \ref{ajustement2} does not modify the fibrations which are already consistent. Afterwards, on these areas, we have definitively ensured the \emph{equalities} $f_i = \phi_{i,j} f_j$. 

For the same reason, at the moment we pass from a stage $\set{a_1 <\cdots <a_k}$ to a stage $\set{a_1 <\cdots <b_{k-1}}$, the fibrations $f_i$ and the diffeomorphisms $\phi_{i,j}$ are definitively fixed on the sets $\Omega_\I$, where $\I$ is a $k$-tuple beginning with an element of $S_{a_1}$. Therefore, on these intersections of order $k$, we have definitively ensured the \emph{equalities} $f_i = \phi_{i,j} f_j$ and all that is done afterwards on intersections of order $k-1$ will not perturb it.

After this adjustment campaign, we have local fibrations $f_i$ on the sets $\Omega_i:=\Omega_i(\kappa/2)$ and diffeomorphisms $\phi_{i,j}$ such that $\phi_{i,j} \circ f_j = f_i$ on $\Omega_i \cap \Omega_j$. The initial estimates still hold, with different constants.

Let us define an equivalence relation: $x$ and $y$ are considered equivalent if there is an index $i$ such that $x$ and $y$ belong to $\Omega_i$ and $f_i(x) = f_i(y)$. Denote by $X$ the quotient topological space and by $\pi$ the corresponding projection. Maps $f_i$ induce homeomorphisms (from their domain to their image) $\check{f_i}$, which endow $X$ with a structure of smooth $3$-manifold: for every (relevant) pair $i,j$, $ \check{f_i} \check{f_j}^{-1} = \phi_{i,j}$ is a diffeomorphism between open sets in $\Rl^{3}$. By construction, $\pi$ is then a smooth fibration.  
\endproof

\subsection{The circle fibration geometry.}
In this whole paragraph, the setting is a complete hyperk\"ahler manifold $(M^4,g)$ with
$$
\int_M \abs{\Rm}^2 r dvol < \infty
$$
and  
$$
\forall \, x \in M, \, \forall \, t \geq 1, \, A t^3 \leq \vol B(x,t) \leq B t^3
$$
($0 < A \leq B$). We have built a circle fibration $\appli{\pi}{M \backslash K}{X}$. The vectors that are tangent to the fibers will be called ``vertical'' whereas vectors orthogonal to the fibers will be called ``horizontal''. Let us average the metric $g$ along the fibers of this fibration. Given a point $x$ in $M \backslash K$, we can choose a unit vector field $V$, defined on a saturated neighborhood of $x$ and vertical (there are two choices of sign). Let $\phi_t$ be the flow of $V$. Denote by $l_x$ the length of the fiber $\pi^{-1}(\pi(x))$. We define a scalar product on $T_x M$ by the formula 
$$
h_x := \frac{1}{l_x} \int_{0}^{l_x} \phi_t^*g \, dt.
$$  
This definition does not depend on the choice of $V$. We thus obtain a Riemannian metric $h$ on $M \backslash K$ and the flows $\phi_t$ are isometries for $h$. To estimate the closeness of $h$ to $g$, we proceed to a few estimations.  

First we show that a local unit vertical field $V$ is almost parallel and almost Killing.

\begin{lem}\label{estimvert}
The covariant derivatives of $V$ can be estimated by 
$$
\abs{\nabla V} \leq C_1 r^{-2} \quad \text{and} \quad
\forall k \geq 2, \, \abs{\nabla^k V} \leq C_k r^{-k}.
$$
\end{lem}

\proof
Let $\appli{f}{\Omega}{\Rl^{3}}$ be one of the local fibrations. By construction, we have $df(V)=0$. Differentiation yields: 
\begin{equation}\label{V1}
\nabla^2 f (V,.) = - df(\nabla V).
\end{equation}
Since $V$ has constant norm, one has 
\begin{equation}\label{N1}
(\nabla V,V)=0
\end{equation}
so, with (\ref{fibrloc}): $\abs{\nabla V} \leq C \abs{\nabla^2 f} \leq C r^{-2}$.
We then make an inductive argument, assuming the result up to order $k-1$. Differentiating $k-1$ times (\ref{V1}), 
we get a formula which looks like
$$
df (\nabla^k V) = \sum_{i=1}^{k-1} \nabla^{1+k-i} f * \nabla^i V 
+ \sum_{i=0}^{k-1} \nabla^{1+k-i} f * \nabla^i V,
$$
which enables us to bound the horizontal part of $\nabla^k V$ by
$$
\abs{\nabla^k V^\bot} \leq C_k \sum_{i=1}^{k-1} \abs{\nabla^{1+k-i} f} \abs{\nabla^i V} 
+ C_k \sum_{i=0}^{k-1} \abs{\nabla^{1+k-i} f} \abs{\nabla^i V}.
$$
Induction assumption and (\ref{fibrloc}) yield: $\abs{\nabla^k V^\bot} \leq C_k (r^{-k} + r^{-k}) \leq C_k r^{-k}$. Differentiating (\ref{N1}), 
we get 
$$
\abs{(\nabla^k V, V)} \leq C_k \sum_{i=1}^{k-1} \abs{\nabla^{k-i} V} \abs{\nabla^i V}, 
$$
so that, by induction assumption: $\abs{(\nabla^k V, V)} \leq C_k r^{-k}$. All in all: $\abs{\nabla^k V} \leq C_k r^{-k}$.
\endproof

\begin{lem}\label{estimlie}
The Lie derivative of $g$ along $V$ satisfies:
$$
\abs{L_V g} \leq C_0 r^{-2},
\quad \text{and} \quad 
\forall \, k \geq 1, \; \abs{\nabla^k L_V g} \leq C_k r^{-1-k}.
$$
\end{lem}

\proof
The formula $L_V g (X,Y) = (\nabla_X V, Y) + (\nabla_Y V,X)$ ensures that for any natural number $k$, $\abs{\nabla^k L_V g}$ is estimated by $\abs{\nabla^{k+1} V}$. So we can apply lemma \ref{estimvert}.
\endproof

If $\phi^t$ is the flow $V$, we are interested in the family of metrics $g_t := \phi^{t*} g$, with Levi-Civita connection $\nabla^t$ and curvature $\Rm^t$. First, a nice formula.

\begin{lem}\label{connexion}
For every vector fields  $X$ and $Y$,  
$$
\frac{d}{dt} \nabla_X^t Y = \Rm^t(X,V)Y - \nabla_{X,Y}^{t, 2} V.
$$
\end{lem}

\proof
The connection $\nabla^t$ is obtained by transporting $\nabla$ thanks to the isometry $\phi^t$:
\begin{equation}\label{formconnexion}
\nabla^t_X Y = \phi^{t *} \nabla_{\phi^t_* X} \phi^t_* Y.
\end{equation}
Thus:
$$
\frac{d}{dt} \phi^t_*\nabla^t_X Y = \frac{d}{dt} \nabla_{\phi^t_* X} \phi^t_* Y,
$$
hence
$$
\phi^t_* [V,\nabla^t_X Y] + \phi^t_* \frac{d}{dt} \nabla^t_X Y = 
\nabla_{[V,\phi^t_* X]} \phi^t_* Y + \nabla_{\phi^t_* X} [V,\phi^t_* Y],
$$
which, thanks to (\ref{formconnexion}) and the invariance of $V$ under its flow, simplifies into 
$$
\frac{d}{dt} \nabla^t_X Y = 
\nabla^t_{[V, X]} Y + \nabla^t_{X} [V, Y]
- [V,\nabla^t_X Y].
$$
Develop and simplify:
$$
\frac{d}{dt} \nabla^t_X Y = 
\nabla^t_{[V, X]} Y + \nabla^t_{X} \nabla^t_V Y - \nabla^t_{X} \nabla^t_Y V
- \nabla^t_V \nabla^t_X Y + \nabla^t_{\nabla^t_X Y} V
= \Rm^t(X,V)Y - \nabla_{X,Y}^{t, 2} V.
$$
\endproof

 This formula gives a control on the covariant derivatives of $g_t$ (with respect to $g$).

\begin{lem}\label{estimflot}
For every $t$, $g_t$ satisfies
$$
\abs{g_t - g} \leq C_0 r^{-2}
\quad \text{and} \quad
\Forall{k}{\Nl^*} \abs{\nabla^k g_t} \leq C_k r^{-1-k}.
$$
\end{lem}

\proof
Let $X$ be a vector field. The definition of the Lie derivative reads:
$$
\frac{d}{dt}g_t(X,X) = (\phi^{t*} L_V g)(X,X).
$$
So, denoting by $\mathcal{L}$ the supremum of $\abs{L_V g}$ on the fiber under consideration, we get
$$
-\mathcal{L} g_t(X,X) \leq \frac{d}{dt} g_t(X,X) \leq \mathcal{L} g_t(X,X).
$$
After integration, we obtain
$$
g(X,X) e^{-\mathcal{L} t} \leq g_t(X,X) \leq g(X,X) e^{\mathcal{L} t}. 
$$
Lemma \ref{estimlie} bounds $\mathcal{L}$:
$$
g(X,X) e^{-C r^{-2}} \leq g_t(X,X) \leq g(X,X) e^{C r^{-2}}, 
$$
hence the first estimate. Now, we consider three vector fields $X$, $Y$, $Z$. We have 
$$
(\nabla^t_X g_t)(Y,Z) = 0 = X \cdot g_t(Y,Z) - g_t(\nabla^t_X Y,Z) -g_t(Y,\nabla^t_X Z),
$$
and
$$
(\nabla_X g_t)(Y,Z) = X \cdot g_t(Y,Z) - g_t(\nabla_X Y,Z) -g_t(Y,\nabla_X Z),
$$
so, if $A^t:=\nabla^t - \nabla$, we arrive at
$$
(\nabla_X g_t)(Y,Z) = g_t(A^t(X,Y),Z) + g_t(Y,A^t(X,Z)),
$$
which we write
\begin{equation}\label{formnablagt}
\nabla g_t = g_t * A^t.
\end{equation}
Lemma \ref{connexion} implies 
$$
A^t = \int_0^t (\Rm^s(.,V) - \nabla^{s,2} V) ds.
$$
Since the curvature is invariant under isometries, we find 
\begin{equation}\label{Rs}
\Rm^t = \phi^{t*} \Rm
\end{equation}
and, thanks to (\ref{formconnexion}) and the invariance of $V$ under the flow, 
\begin{equation}\label{nabla2Vs}
\nabla^{t,2} V = \phi^{t*} \nabla^{2} V .
\end{equation}
We have bounds on $g_t$, $\Rm$ and $\nabla^{2} V$: (\ref{formconnexion}) leads to $\abs{\Rm^t} \leq C r^{-2}$ (and even $r^{-3}$), $\abs{\nabla^{t,2} V} \leq C r^{-2}$ and $\abs{A^t} \leq C r^{-2}$. 

Now assume (by induction) that for some $k\geq 1$, for every $0 \leq i \leq k-1$ and every $t$,
\begin{eqnarray*}
\abs{\nabla^i (g_t - g)} &\leq& C r^{-1-i}, \\
\abs{\nabla^i \Rm^t} &\leq& C r^{-2-i}, \\
\abs{\nabla^i \nabla^{t,2} V} &\leq& C r^{-2-i}.
\end{eqnarray*}
In particular, we get 
$$
\forall \, t, \, \Forall{i}{[0,k-1]} \abs{\nabla^i A^t} \leq C r^{-2-i}.
$$
Fix $t$. Differentiating (\ref{formnablagt}), we obtain the formula
$$
\nabla^k g_t = \sum_{i=0}^{k-1} \nabla^{k-1-i} g_t * \nabla^i A^t.
$$
Induction assumption yields $\abs{\nabla^k g_t} \leq C r^{-1-k}$. To go on, we need to estimate $\abs{\nabla^{t,i} A^t}$, $i \leq k-1$. To do this, we write $\nabla^t = \nabla + A^t$. In this way, we see that $\abs{\nabla^{t,i} A^t}$ can be controlled by a sum of a bounded number of terms like 
$$
\left( \prod_{\alpha=0}^{i-1} \abs{\nabla^{\alpha} A^t}^{m_\alpha} \right) \abs{\nabla^{\beta}A^t}
$$
with natural numbers $m_\alpha$, $\beta$ satisfying
$$
\sum_{\alpha=0}^{i-1} (1+\alpha)m_\alpha + \beta = i.
$$
Induction assumption implies each of these terms is bounded by 
$C r^{-(2+\alpha)m_\alpha -2-\beta}\leq C r^{-2-i}$, so 
$$
\abs{\nabla^{t,i} A^t} \leq C r^{-2-i}.
$$
Now, writing $\nabla = \nabla^t - A^t$, we estimate $\abs{\nabla^k \Rm^t}$ by a sum of a bounded number of terms like
$$
\left( \prod_{\alpha=0}^{k-1} \abs{\nabla^{t,\alpha} A^t}^{m_\alpha} \right) \abs{\nabla^{t,\beta}\Rm^t}
$$
with natural numbers $m_\alpha$, $\beta$ satisfying
$$
\sum_{\alpha=0}^{k-1} (1+\alpha)m_\alpha + \beta = k.
$$
With (\ref{Rs}) and (\ref{formconnexion}), we bound $\abs{\nabla^{t,\beta}\Rm^t}$ by $\abs{\nabla^{\beta}\Rm}$ and thus by $r^{-2 -\beta}$. Eventually, we find $\abs{\nabla^k \Rm^t} \leq C r^{-2-k}$. In the same way, we get $\abs{\nabla^k \nabla^{t,2} V} \leq C r^{-2-k}$ and we conclude by induction.
\endproof

\begin{lem}\label{estimlong}
The length $l$ of the fibers is controlled by:
$$
\abs{dl} \leq C_1 r^{-2} \quad \text{and} 
\quad \forall k\geq 2, \, \abs{\nabla^k l} \leq C_k r^{-k}.
$$
\end{lem}

\proof
By construction, we have the identity $\phi^{l(x)}(x) = x$,
at every point $x$ in $M\backslash K$. Differentiation yields $dl \otimes V + T \phi^l = id$. Taking the scalar product 
with $V$, we obtain $dl = (g- g_l)(V,.)$. Differentiating this leads to
$$
\nabla^k l = \sum_{i=0}^{k-1} \nabla^i (g-g_l) * \nabla^{k-1-i}V.
$$
Now we use (\ref{estimflot}), (\ref{estimvert}) and the bound on $l$: $\abs{\nabla^k l} \leq C r^{-k}$.
\endproof

We can finally control the metric $h$, obtained by averaging $g$ along the fibers.

\begin{prop}\label{estimh}
The averaged metric $h$ obeys the estimates:
$$
\abs{h-g} \leq C_k r^{-2}
\quad \text{and} \quad
\forall k \geq 1, \, \abs{\nabla^k h} \leq C_k r^{-1-k}.
$$
\end{prop}

\proof
The definition of $h$ can be written
$$
h-g = \frac{1}{l} \int_0^l (g_t-g) dt 
$$
The first estimate follows immediately from (\ref{estimflot}). Let us differentiate:
$$
\nabla h = \frac{dl}{l} \otimes (g_l - h) + \frac{1}{l} \int_0^l \nabla g_t dt. 
$$
An induction yields for every $k \geq 1$:
$$
\nabla^k h = \sum_{i=1}^{k} C_{k}^i \frac{\nabla^i l}{l} \otimes \nabla^{k-i}(g_l-h) 
+ \frac{1}{l} \int_0^l \nabla^k g_t dt. 
$$
(\ref{estimflot}) and (\ref{estimlong}) then lead, by induction, to $\abs{\nabla^k h} \leq C r^{-1-k}$. 
\endproof

Since $g$ has cubic curvature decay, we deduce the

\begin{cor}\label{courbh}
The curvature of $h$ has cubic decay.
\end{cor}

Now, let us push $h$ down into a Riemannian metric $\check{h}$ on $X$: for every point $y$ in $X$, for every vector $w$ in $T_y X$, we choose a lift $x$ of $y$ ($\pi(x) =y$) and we set $\check{h}_y(w,w) = h_x(v,v)$ where $v$ is the horizontal lift of $w$ in $T_x M$; this definition makes sense because the flow $\phi_t$ is isometric for $h$.    

\begin{prop}
The manifold $X$ is diffeomorphic to the complementary set of a ball in $\Rl^{3}$, mod out by the action of a finite subgroup of $O(3)$. Moreover, $\check{h}$ is an ALE metric of order $1^-$, that is
$$
\check{h}= g_{\Rl^3} + \OO(r^{-\tau}) \quad \text{for every } \tau<1.
$$
\end{prop}

\proof
Observe the volume of a ball of radius $t$ in $(X^3,\check{h})$ is comparable to $t^{3}$. To estimate the curvature on the base, we use O'Neill formula (\cite{Bes}), which asserts that if $Y$ and $Z$ are orthogonal unit horizontal vector fields on $M \backslash K$, then 
$$
\Sect_{\check{h}}(\pi_* Y \wedge \pi_* Z) = \Sect_{h}(Y \wedge Z) + \frac{3}{4} h([Y,Z],V)^2. 
$$
The first term decays at a cubic rate by \ref{courbh}. Moreover, 
$$
h([Y,Z],V) = - (\nabla_Y h)(Z,V) - h(Z,\nabla_Y V) + (\nabla_Z h)(Y,V) + h(Y,\nabla_Z V). 
$$ 
Lemma \ref{estimvert} and corollary \ref{estimh} yield
$
\abs{h([Y,Z],V)} \leq C r^{-2}.
$
Hence:
$$
\abs{\Sect_{\check{h}}(\pi_* Y \wedge \pi_* Z)} \leq C r^{-3}.
$$
This cubic curvature decay, combined with Euclidean volume growth, enables us to apply the main theorem of \cite{BKN}.
\endproof


\subsection{What have we proved ?}

We have proved the following theorem.

\begin{thm}\label{thinstanton}
Let $(M^4,g)$ be a complete hyperk\"ahler manifold satisfying
$$
\int_M \abs{\Rm}^2 r \, dvol < \infty
$$
and
$$
\forall \, x \in M, \, \forall \, t \geq 1, \, A t^\nu \leq \vol B(x,t) \leq B t^\nu
$$
with $0 < A \leq B$ and $3 \leq \nu <4$. Then there is compact set $K$ in $M$, a ball $B$ in $\Rl^{3}$, a finite subgroup $G$ of $O(3)$ and a circle fibration $\appli{\pi}{M \backslash K}{(\Rl^{3} \backslash B)/ G}$. Moreover, the metric $g$ obeys
$$
g = \pi^*\tilde{g} + \eta^2 + \OO(r^{-2}),
$$  
where $\eta^2$ measures the projection along fibers and $\tilde{g}$ is an ALE metric of order $1^-$.
\end{thm}

Let us precise the topology at infinity, that is the topology of the connected space $E=M \backslash K$, which, according to theorem \ref{thinstanton}, is a circle bundle over $X = \Rl^{3}\backslash B /G$. Thanks to the projection $\appli{p}{\bar{X} = \Rl^{3}\backslash B}{X = \Rl^{3}\backslash B /G}$, we can pull back the fibration $\pi$ into a circle fibration $\appli{\bar{\pi}}{\bar{E}}{\bar{X}}$. The space $\bar{E}$ is a finite covering of $E$, with order $\abs{G}$:
$$
\bar{E} = \set{(\bar{x},e) \in \bar{X} \times E, \, p(\bar{x}) = \pi(e)}
$$ 
and $\bar{\pi}$ is given by the projection onto the first factor ($pr_1$). 
\diagcommut{\bar{E}}{pr_2}{E}{\bar{\pi}}{\pi}{\bar{X}}{p}{X}
  
Of course, $\bar{X} = \Rl^{3}\backslash B$ has the homotopy type of $\Sph^{2}$, so that we can classify its circle fibrations. Moreover, the homotopy groups of $\bar{E}$ can be computed thanks to the long exact homotopy sequence associated to $\bar{\pi}$. In this way, we obtain essentially two cases, which are distinguished by the homotopy groups at infinity (those of $M \backslash K$).

\begin{itemize}
\item If the fundamental group at infinity is finite, then a finite covering of $M \backslash K$ is
$\Rl^4 \backslash \mathbb{B}^4$ and the circle fibration is the Hopf fibration, up to a finite group action. In this case, the $\pi_2$ at infinity is trivial. This is typically the ``Taub-NUT'' situation .
\item If the fundamental group at infinity is infinite, then a finite covering of $M \backslash K$ is
$\Rl^3 \backslash \mathbb{B}^3 \times \Sph^1$ and the circle fibration comes from the trivial one. The $\pi_2$ at infinity is then $\Ir$. 
\end{itemize}

It is easy to adapt the arguments above in order to obtain the following result.

\begin{thm}\label{thmgl}
Let $(M^n,g)$ be a complete manifold satisfying
$$
\forall k \in \Nl, \, \abs{\nabla^k \Rm} = \OO(r^{-3-k}) 
$$
and
$$
\forall \, x \in M, \, \forall \, t \geq 1, \, A t^{n-1} \leq \vol B(x,t) \leq \omega(t) t^n
$$
for some positive number $A$ and some function $\omega$ going to zero at infinity. Assume 
moreover there is $c \geq 1$ such that the holonomy $H$ of any geodesic loop based at $x$ and with length $L \leq r(x)/c$ satisfies
$$
\abs{H - \id} \leq \frac{c L}{r(x)}.
$$  
Then there is compact set $K$ in $M$, a ball $B$ in $\Rl^{n-1}$, a finite subgroup $G$ of $O(n-1)$ and a circle fibration $\appli{\pi}{M \backslash K}{(\Rl^{n-1} \backslash B)/ G}$. Moreover, the metric $g$ obeys
$$
g = \pi^*\tilde{g} + \eta^2 + \OO(r^{-2}),
$$  
where $\eta^2$ measures the projection along fibers and $\tilde{g}$ is an ALE metric of order $1-$ ($1$ if $n \geq 5)$. 
\end{thm}

\begin{rem}
The required estimates on the curvature are satisfied on a Ricci flat manifold with cubic curvature decay. This allows one to englobe the Schwarzschild metrics (\cite{Min} for instance) 
in this setting. Note a little topology ensures the fibration is trivial if $n \geq 5$.  
\end{rem}

\appendix

\section{Curvature decay.}

The following result is proved in \cite{Min}. Recall we always distinguish a point $o$ in our manifolds. We will use the
measure $\mu$ defined by $d\mu = \frac{r^n}{\vol B(o,r)} dvol$. 

\begin{thm}
Let $(M^n,g)$ be a complete Ricci flat manifold. Assume there are numbers $\nu > 2$ and  $C >0$ such that 
$$
\forall t \geq s > 0, \  \frac{\vol B(o,t)}{\vol B(o,s)} \geq C \left( \frac{t}{s} \right)^\nu.
$$
Then the integral bound
$$
\int_M \abs{\Rm}^{\frac{n}{2}} d\mu < \infty
$$
implies the pointwise bound
$$
\abs{\Rm} = \OO(r^{-a(n,\nu)}) \quad \text{ with }\quad   a(n,\nu) = \max \left( 2 , \frac{(\nu-2)(n-1)}{n-3} \right).
$$ 
\end{thm}

It should be stressed that the integral assumption 
$\int \abs{\Rm}^{\frac{n}{2}} d\mu < \infty$ is weaker than $\abs{\Rm} = \OO(r^{-2-\epsilon})$ for some positive $\epsilon$.
In this paper, we use the

\begin{cor}
Let $(M^n,g)$ be a complete Ricci flat manifold, with $n \geq 4$. Assume there are positive numbers $A$ and  $B$ such that 
$$
\forall t \geq 1, \  A t^{n-1} \leq \vol B(o,t) \leq B t^{n-1}.
$$
Then the integral bound
$$
\int_M \abs{\Rm}^{\frac{n}{2}} r \, dvol < \infty
$$
implies the pointwise bound
$$
\abs{\Rm} = \OO(r^{-(n-1)}).
$$ 
\end{cor}

These estimates follow from a Moser iteration, which is self improved thanks to a global weighted Sobolev inequality. In this appendix,
we wish to obtain similar estimates on the covariant derivatives of the curvature tensor. To do this, we need a technical 
inequality. 

\begin{lem}[Moser iteration with source term]\label{moser}
Let $(M^n,g)$ be a complete noncompact Riemannian manifold with nonnegative Ricci curvature and let $E \To M$ be a smooth Euclidean vector bundle, endowed with a compatible connection $\nabla$. We denote by $\overline{\Delta} = \nabla^* \nabla$ the Bochner Laplacian and suppose $V$ is a continuous field 
of symmetric endomorphisms of $E$ whose negative part satisfies $\abs{V_-} = \OO (r^{-2})$. Given a locally bounded section $\phi$ and a locally Lipschitz section $\sigma$ such that $(\sigma,\overline{\Delta} \sigma + V \sigma) \leq (\sigma,\phi)$, the following estimate holds for large $R$:
$$
\sup_{A(R,2R)} \abs{\sigma} 
\leq \frac{C  }{\vol B(o,R)^{\frac{1}{2}}}  \norm{\sigma}_{L^2(A(R/2,5R/2))} + C R^2  \norm{\phi}_{L^\infty(A(R,2R))}.
$$  
\end{lem}

\proof
Set $u := \abs{\sigma} + F$, with 
$
F := R^2 \norm{\phi}_{L^\infty(A(R,2R))}.
$
The case $\phi = 0$ is treated in \cite{Min}. Actually, in \cite{Min}, the estimation is written assuming a global weighted Sobolev inequality. But since we work at a fixed scale $R$, there is no need for such a global inequality: the local Sobolev inequality of L. Saloff-Coste \cite{SC}, with controlled constant, is sufficient for our purpose; and its validity only requires $\Ric \geq 0$. Therefore we assume $F \not= 0$. 

To avoid troubles on the zero set of $\sigma$, let us consider the regularizations $v_\epsilon := \sqrt{\abs{\sigma}^2 + \epsilon}$ and 
$u_\epsilon := v_\epsilon + F$. Observing the inequalities
$$
v_\epsilon \Delta v_\epsilon \leq (\sigma,\overline{\Delta} \sigma) \leq \abs{\sigma} ( \abs{V_-} \abs{\sigma} + \abs{\phi} \abs{\sigma})
\leq v_\epsilon ( \abs{V_-} \abs{\sigma} + \abs{\phi} \abs{\sigma}),
$$
we deduce
$
\Delta v_\epsilon \leq \abs{V_-} v_\epsilon + \abs{\phi}
$
and thus find 
$$
\Delta u_\epsilon \leq \abs{V_-} u_\epsilon + \abs{\phi} \leq \left(\abs{V_-}  + \frac{\abs{\phi}}{F} \right) u_\epsilon.
$$
Our choice of $F$ enables us to use the estimate without source term in \cite{Min}:
$$
\sup_{A(R,2R)} u_\epsilon 
\leq \frac{C  }{\vol B(o,R)^{\frac{1}{m}}} \norm{u_\epsilon}_{L^m(A(R/2,5R/2))}
$$
Let $\epsilon$ go to zero, so as to obtain
$$
\sup_{A(R,2R)} \abs{\sigma} \leq \sup_{A(R,2R)} u 
\leq \frac{C  }{\vol B(o,R)^{\frac{1}{m}}} \norm{\sigma}_{L^m(A(R/2,5R/2))} + C F,
$$
which is what we want.
\endproof

We will use this lemma on tensor bundles, with the induced Levi-Civita connection, in order to prove that on a Ricci flat manifold, if the curvature decays at infinity, then the covariant derivatives of the curvature also decay.

\begin{prop}
Let $(M^n,g)$ be a complete noncompact Ricci flat manifold. If $a\geq 2$, the estimate 
$
\abs{\Rm} = \OO (r^{-a})
$ 
implies for positive integer $i$:
$
\abs{\nabla^i \Rm} = \OO(r^{-a-i}).
$ 
\end{prop}

\proof
Since $M$ is Ricci flat, its curvature tensor obeys an elliptic equation \cite{BKN}
$$
\overline{\Delta} \Rm = \Rm * \Rm,
$$
which implies for every $k$ in $\Nl$ \cite{TV}:
\begin{equation}\label{formule}
\overline{\Delta} \nabla^{k} \Rm = \sum_{i=0}^{k} \nabla^i \Rm * \nabla^{k-i} \Rm.
\end{equation}
Let us prove the result by induction on $i$. The case $i=0$ is contained in the assumptions. Suppose that the result is established for $i \leq k$, with $k \geq 0$. Formula (\ref{formule}) can be written  
$$
(\overline{\Delta} - \Rm *) \nabla^{k+1} \Rm =  \sum_{i=1}^{k} \nabla^i \Rm * \nabla^{k+1-i} \Rm.
$$ 
Since the right-hand side is bounded by $C_{k+1} r^{-2a-k-1}$, lemma \ref{moser} yields:  
\begin{equation}\label{estimation}
\sup_{A(R,2R)} \abs{\nabla^{k+1} \Rm} 
\leq \frac{C_{k+1}}{\vol B(o,R)^{\frac{1}{2}}}  \norm{\nabla^{k+1} \Rm}_{L^2(A(R/2,5R/2))} + C_{k+1} R^{1-2a-k}.
\end{equation}
Let $\chi$ be a positive smooth function equal to $1$ on $A(R/2,5R/2)$, $0$ on $A(R/3,3R)^c$ and with differential bounded by $10/R$. Then 
we can write 
$$
\int_{A(R/2,5R/2)} \abs{\nabla^{k+1} \Rm}^2 \leq \int_{A(R/3,3R)} \abs{\nabla \left(\chi\nabla^{k} \Rm \right)}^2   
$$
and, after integration by parts, we find
$$
\int_{A(R/2,5R/2)} \abs{\nabla^{k+1} \Rm}^2  \leq \int_{A(R/3,3R)} \abs{d\chi}^2 \abs{\nabla^{k} \Rm }^2  
+ \int_{A(R/3,3R)} \chi^2 (\nabla^{k} \Rm,\overline{\Delta} \nabla^{k} \Rm) .  
$$
With (\ref{formule}), we obtain the upper bound
\begin{eqnarray*}
\int_{A(R/2,5R/2)} \abs{\nabla^{k+1} \Rm}^2  &\leq& \frac{100}{R^2}  \int_{A(R/3,3R)} \abs{\nabla^{k} \Rm }^2  \\
&+& C_{k+1} \sum_{i=0}^{k} \int_{A(R/3,3R)} \abs{\nabla^{k} \Rm} \abs{\nabla^{i} \Rm} \abs{\nabla^{k-i} \Rm}.   
\end{eqnarray*}
Using $a\geq 2$, we estimate this by
\begin{eqnarray*}
\int_{A(R/2,5R/2)} \abs{\nabla^{k+1} \Rm}^2  &\leq& C_{k+1} \vol B(o,R) \left( R^{-2-2a-2k}
+ R^{-3a-2k}  \right) \\
&\leq& C_{k+1} \vol B(o,R) R^{-2-2a-2k}.
\end{eqnarray*}
As a result, (\ref{estimation}) implies
$$
\sup_{A(R/2,5R/2)} \abs{\nabla^{k+1} \Rm} \leq C_{k+1} \left( R^{-1-a-k} +  R^{1-2a-k} \right) \leq C_{k+1}  R^{-1-a-k},
$$
hence
$$
\abs{\nabla^{k+1} \Rm} \leq C_{k+1}  r^{-a-(k+1)}.
$$
\endproof

\begin{cor}
Let $(M^n,g)$ be a complete Ricci flat manifold. Assume there are numbers $\nu > 2$ and  $C >0$ such that 
$$
\forall t \geq s > 0, \  \frac{\vol B(o,t)}{\vol B(o,s)} \geq C \left( \frac{t}{s} \right)^\nu.
$$
Then the integral bound
$$
\int_M \abs{\Rm}^{\frac{n}{2}} d\mu < \infty
$$
implies for every $k$ in $\Nl$:
$$
\abs{\nabla^k\Rm} = \OO(r^{-a(n,\nu)-k}) \quad \text{ with }\quad   a(n,\nu) = \max \left( 2 , \frac{(\nu-2)(n-1)}{n-3} \right).
$$ 
\end{cor}

\begin{cor}\label{curvdec}
Let $(M^n,g)$ be a complete Ricci flat manifold, with $n \geq 4$. Assume there are positive numbers $A$ and  $B$ such that 
$$
\forall t \geq 1, \  A t^{n-1} \leq \vol B(o,t) \leq B t^{n-1}.
$$
Then the integral bound
$$
\int_M \abs{\Rm}^{\frac{n}{2}} r \, dvol < \infty
$$
implies for every $k$ in $\Nl$:
$$
\abs{\nabla^k\Rm} = \OO(r^{-(n-1)-k}).
$$ 
\end{cor}


\section{Distance and curvature.}

The following lemma sums up some comparison estimates on the distance function. Up to order two, it is quite classical. Higher order estimates do not seem to be proved in the litterature, so we include a proof.  

\begin{lem}\label{derdist}
Consider a complete Riemannian manifold $(M,g)$, a point $x$ in $M$ and a number $a \geq 2$ such that 
$$
\inj(x) > \epsilon \geq 1
$$ 
and
$$
\Forall{i}{[0,k]} \abs{\nabla^i \Rm} \leq c \epsilon^{-a-i}    
$$
on the ball $B(x,\epsilon)$. Then there is a constant $C$ such that on this ball, the function $\rho=d(x,.)^2/2$ satisfies:
\begin{itemize}
\item $\abs{d\rho} \leq \epsilon$ ;
\item $\abs{\nabla^2 \rho - g} \leq C \epsilon^{2-a}$ ;
\item for $3 \leq i \leq k$, $\abs{\nabla^i \rho} \leq C \epsilon^{4-a-i}$. 
\end{itemize}
\end{lem}

\proof
The first estimate is obvious and the second follows from \cite{BK}. Let us turn to higher order estimates. We consider the gradient $N$ of 
$r := d(x,.)$ and use the Riccati equation for the second fundamental form $\nabla N$ of geodesic spheres:
$$
\nabla_N S = - S^2 - \Rm(N,.)N.
$$
Identifying quadratic forms to symmetric endomorphisms, we write the endomorphism $E :=\nabla^2 \rho - \Id$ as 
$$
E = dr \otimes N + r S - \Id.
$$
Setting $V=\grad \rho = r N$, we obtain the equation 
$$
\nabla_V E = - E - E^2 - \Rm(V,.)V.
$$
Since $\nabla V = \Id + E$ and 
$$
\nabla_V \nabla E = \nabla \nabla_V E - \nabla_{\nabla V} E + \Rm(V,.) E,
$$
it follows that
$$
\nabla_V \nabla E = -2 \nabla E + E * \nabla E + \nabla \Rm * V * V + \Rm * \nabla V * V + \Rm * V.
$$
Observing that for $k \geq 2$, we have $\nabla^kV = \nabla^{k-1} E$, an induction yields: 
\begin{eqnarray*}
\nabla_V \nabla^k E = &-& (k+1) \nabla^k E + \sum_{i+j=k} \nabla^i E * \nabla^j E \\
&+& \sum_{i+j+l=k} \nabla^i \Rm * \nabla^j V * \nabla^l V + \sum_{i+j=k-1} \nabla^i \Rm * \nabla^j V,
\end{eqnarray*}
for every natural number $k$. We set $F_k = r^{k+1} \nabla^k E$ and $G = E/r$, so that
$$
\nabla_N F_k = G * F_k + H_k
$$
with 
$$
H_k = r^{-2} \sum_{i=1}^{k-1} F_i * F_{k-i} 
+ r^k \left( \sum_{i+j+l=k} \nabla^i \Rm * \nabla^{j+1} \rho * \nabla^{l+1} \rho + \sum_{i+j=k-1} \nabla^i \Rm * \nabla^{j+1} \rho
\right)
$$
Along a geodesic starting from $x$, we find
$$
\partial_r \abs{F_k} \leq C_k \abs{F_k} \abs{G} + \abs{H_k}
$$  
and since the order two estimate ensures $r \abs{G}$ is small, we can bound $\abs{F_k}$ by $r \sup \abs{H_k}$, up to a constant. We will prove by induction the estimate
$$
\abs{F_k} \leq C_k r^{k+1} \epsilon^{2-a-k}
$$   
that is 
$$
\abs{\nabla^k E} \leq C_k \epsilon^{2-a-k},
$$   
or  
$$
\abs{\nabla^{k+2} \rho} \leq C_k \epsilon^{4-a-(k+2)}.
$$
It will conclude the proof. Initialization ($k=0$) follows from the order two estimate on $\rho$. Assume the estimates up to order $k-1$. It follows that $H_k$, up to some constant, is bounded by:
$$
r^{-2} r^{k+2} \epsilon^{4-2a-k}  
+ r^k \left( \epsilon^{-a-k+4-a-1+4-a-1} + \epsilon^{-a-k+1+4-a-1} \right)
$$
hence 
$$
\abs{H_k} \leq C_k r^k \left(\epsilon^{4-2a-k} + \epsilon^{4-2a-k+2-a} + \epsilon^{4-2a-k} \right).
$$
With $a \geq 2$, we find $\abs{H_k} \leq C_k r^k \epsilon^{4-2a-k} \leq C_k r^k \epsilon^{2-a-k}$ and therefore we get the promised estimate $\abs{F_k} \leq C_k r^{k+1} \epsilon^{2-a-k}$, hence the result.
\endproof

H. Kaul \cite{Kau} proved a control on Christoffel coefficients in the exponential chart, given bounds on $\Rm$ and $\nabla \Rm$. We need the following 

\begin{prop}\label{kaul}
Consider a complete Riemannian manifold $(M,g)$, a point $x$ in $M$ and a number $a \geq 2$ such that 
$$
\abs{\Rm} \leq c \epsilon^{-a} \quad \text{ and } \abs{\nabla \Rm}  \leq c \epsilon^{-a-1} 
$$
on the ball $B(x,\epsilon)$, with $\epsilon \geq 1$. Then there is a constant $C$ such that, on the ball $\hat{B}(0,\epsilon)$ in $T_x M$, 
the connection $\nabla^{\hat{g}}$ of the metric $\hat{g}=\exp_x^* g$ and the flat connection $\nabla^0$ are related by
$$
\abs{\nabla^{\hat{g}} - \nabla^0} \leq C \epsilon^{1-a}.
$$ 
\end{prop}

A better control on the distance function stems from this.

\begin{lem}\label{compderdist}
Consider a complete Riemannian manifold $(M,g)$, a point $x$ in $M$ and a number $a \geq 2$ such that 
$$
\abs{\Rm} \leq c \epsilon^{-a} \quad \text{ and } \abs{\nabla \Rm}  \leq c \epsilon^{-a-1} 
$$
on the ball $B(x,\epsilon)$, with $\epsilon \geq 1$. Then there is a constant $C$ such that if $v$ and $w$ belong to $\hat{B}(0, C^{-1} \epsilon)$, endowed with $\hat{g}$, then
$$
\abs{(d\rho_v)_w - g_x(w-v,.)} \leq C \epsilon^{3-a}.
$$
\end{lem}

\proof
First, choose a sufficiently large $C$ to ensure the convexity of the ball under consideration. Observe the expression 
$$
(d\rho_v)_w = - \hat{g}_w (\Exp_w^{-1} v, .),
$$
where $\Exp$ is the exponential map of $\hat{g}$. Comparison yields 
\begin{equation}\label{compmet}
\abs{\hat{g}_w - g_x} \leq C \epsilon^{-a} \epsilon^2 = C \epsilon^{2-a}.
\end{equation}
Suppose $\gamma$ parameterizes the geodesic connecting $w$ to $v$ in unit time. The geodesic equation $\nabla^{\hat{g}}_{\dot{\gamma}} \dot{\gamma} = 0$ can be written
$$
\ddot{\gamma} + \left( \nabla^{\hat{g}}_{\dot{\gamma}} - \nabla^0 \right) \dot{\gamma} = 0.
$$
With (\ref{kaul}), we obtain
$$
\abs{\ddot{\gamma}} \leq C \epsilon^{1-a} \epsilon^2 = C \epsilon^{3-a}.
$$
Taylor formula
$$
\gamma(1) - \gamma(0) - \dot{\gamma}(0) = \int_0^1 (1-t) \ddot{\gamma}(t) dt
$$
then yields 
$$
\abs{v - w - \Exp_w^{-1} v} \leq C \epsilon^{3-a}.
$$
To conclude, we write
\begin{eqnarray*}
\abs{(d\rho_v)_w - g_x(w-v,.)} & = & \abs{\hat{g}_w (\Exp_w^{-1} v, .) - g_x(v-w,.)} \\
& \leq & \abs{(\hat{g}_w - g_x) (\Exp_w^{-1} v,.)} 
 +  \abs{g_x (\Exp_w^{-1} v, .) - g_x(v-w,.)} \\
&\leq& C \epsilon^{2-a} \epsilon + C \epsilon^{3-a} \\
&\leq& C \epsilon^{3-a}.
\end{eqnarray*}
\endproof



\end{document}